\RequirePackage{fix-cm}
\documentclass[smallextended]{svjour3}       
\smartqed  
\usepackage{graphicx}

 \usepackage{mathptmx}      
\usepackage{latexsym}

\usepackage[USenglish]{babel}
\usepackage[T1]{fontenc}
\usepackage[ansinew]{inputenc}

\usepackage{lmodern} 

\usepackage{graphicx}
\usepackage{wrapfig}
\usepackage{caption}
\usepackage{subfig}
\usepackage{wrapfig}
\usepackage{pgfplots}

\usepackage{amssymb, amsmath}
\usepackage{amsfonts}
\usepackage{dsfont}

\usepackage{authblk}

\usepackage{calc}
\usepackage{xcolor}
\usepackage{tikz}
\usetikzlibrary{decorations.markings}
\usetikzlibrary{arrows}
\tikzstyle{vertex}=[circle, draw, inner sep=0pt, minimum size=6pt]

\tikzstyle{vertexr}=[fill=red, circle, draw, inner sep=0pt, minimum size=6pt]

\tikzstyle{vertexb}=[fill=blue, circle, draw, inner sep=0pt, minimum size=6pt]

\tikzstyle{vertexg}=[fill=green, circle, draw, inner sep=0pt, minimum size=6pt]

\tikzstyle{vertexy}=[fill=yellow, circle, draw, inner sep=0pt, minimum size=6pt]

\usepackage{color}

\definecolor{gray50}{gray}{0.5}

\usepackage[colorlinks=true]{hyper ref}
\usepackage{enumitem}

\usepackage{multirow}

\usepackage{authblk}

\newcommand{\cible}{\overline{\vec{x}}}

\newcommand{\Pro}{\mathbb{P}}
\newcommand{\E}{\mathbb{E}}

\newcommand{\diag}{{\rm diag}}

\spnewtheorem{notation}{Notation}{\it}{\rm}
\spnewtheorem{assumption}{Hypothesis}{\it}{\rm}
\spnewtheorem{modassumption}{Modeling assumption}{\it}{\rm}
\spnewtheorem{warning}{Warning}{\it}{\rm}
\spnewtheorem*{proposition*}{Proposition}{\bf}{\it}

\smartqed

\begin{document}

\title{Multi-type Galton-Watson processes with affinity-dependent selection applied to antibody affinity maturation
}

\titlerunning{Multi-type Galton-Watson processes with affinity-dependent selection}        

\author{Irene Balelli  \and Vuk Mili\v{s}i\'c    \and Gilles Wainrib
}


\institute{I. Balelli  \and V. Mili\v{s}i\'c \at
              Universit\'e Paris 13, Sorbonne Paris Cit\'e, LAGA, CNRS (UMR 7539). \\
F-93430 - Villetaneuse - France. 
              \email{milisic@math.univ-paris13.fr }           
	\and 
	I. Balelli \at ISPED, Centre INSERM U1219, and INRIA - Statistics in System Biology and Translational Medicine Team. F-33000 - Bordeaux - France. 
	\email{irene.balelli@inserm.fr}
           \and
           G. Wainrib \at
              Ecole Normale Sup\'erieure, D\'epartement d'Informatique. \\
              45 rue d'Ulm, 75005 - Paris - France. 
              \email{gilles.wainrib@ens.fr}
}

\date{Received: date / Accepted: date}

\maketitle

\begin{abstract}

We analyze the interactions between division, mutation and selection in a simplified evolutionary model, assuming that the population observed can be classified into fitness levels. The construction of our mathematical framework is motivated by the modeling of antibody affinity maturation of B-cells in Germinal Centers during an immune response. This is a key process in adaptive immunity leading to the production of high affinity antibodies against a presented antigen. Our aim is to understand how the different biological parameters affect the system's functionality. We identify the existence of an optimal value of the selection rate, able to maximize the number of selected B-cells for a given generation.



\end{abstract}

\tableofcontents

\section{Introduction}\label{sec:intro}

Antibody Affinity Maturation (AAM) takes place in Germinal Centers (GCs), specialized micro-environnements which form in the peripheral lymphoid organs upon infection or immunization \cite{victora2014clonal,de2015dynamics}. GCs are seeded by ten to hundreds distinct B-cells \cite{tas2016visualizing}, activated after the encounter with an antigen, which initially undergo a phase of intense proliferation \cite{de2015dynamics}. Then, AAM is achieved thanks to multiple rounds of division, Somatic Hypermutation (SHM) of the B-cell receptor proteins, and subsequent selection of B-cells with improved ability of antigen-binding \cite{maclennan2000b}. B-cells which successfully complete the GC reaction output as memory B-cells or plasma cells \cite{victora2012germinal,de2015dynamics}. Indirect evidence suggests that only B-cells exceeding a certain threshold of antigen-affinity differentiate into plasma cells \cite{phan2006high}. The efficiency of GCs is assured by the contribution of other immune molecules, for instance Follicular Dendritic Cells (FDCs) and follicular helper T-cells (Tfh). Nowadays the key dynamics of GCs are well characterized \cite{maclennan2000b,de2015dynamics,gitlin2014clonal,tas2016visualizing}. Despite this there are still mechanisms which remain unclear, such as the dynamics of clonal competition of B-cells, hence how the selection acts. In recent years a number of mathematical models of the GC reaction has appeared to investigate these questions, such as \cite{meyer2006analysis,wang2016stochastic}, where agent-based models are developed and analyzed through extensive numerical simulations, or \cite{zhang2010optimality} where the authors establish a coarse-grained model, looking for optimal values of \emph{e.g.} the selection strength and the initial B-cell fitness maximizing the affinity improvement.\\

Our aim in this paper is to contribute to the mathematical foundations of adaptive immunity by introducing and studying a simplified evolutionary model inspired by AAM, including division, mutation, affinity-dependent selection and death. We focus on interactions between these  mechanisms,  identify and analyze the parameters which mostly influence the system functionality, through a rigorous mathematical analysis. 
This research is motivated by important biotechnological applications. 
Indeed, the fundamental understanding of the evolutionary mechanisms involved in AAM have been inspiring many methods for the synthetic production of specific antibodies for drugs, vaccines or cancer immunotherapy \cite{ansari2010identification,kringelum2012reliable,shen2012towards}. This production process involves the selection of high affinity peptides and requires smart methods to ge\-ne\-rate an appropriate diversity \cite{currin2015synthetic}. Beyond biomedical motivations, the study of this learning process has also given rise in recent years to a new class of bio-inspired algorithms \cite{LNDeCaFJVoZu,pang2015clonalal,timmis2008theoretical}, mainly addressed to solve optimization and learning problems.\\

We consider a model in which B-cells are classified into $N+1$ affinity classes with respect to a presented antigen, $N$ being an integer big enough to opportunely describe the possible fitness levels of a B-cell with respect to a specific antigen \cite{weiser2011affinity,xu2015key}. A B-cell is able to increase its fitness thanks to SHMs of its receptors: only about $20\%$ of all mutations are estimated to be affinity-affecting mutations \cite{MShaRMeh,shlomchik1998clone}. By conveniently define a transition probability matrix, we can characterize the probability that a B-cell belonging to a given affinity class passes to another one by mutating its receptors thanks to SHMs. Therefore we define a selection mechanism which acts on B-cells differently depending on their fitness. We mainly focus on a model of positive and negative selection in which B-cells submitted to selection either die or exit the GC as output cells, according to the strength of their affinity with the antigen. Hence, in this case, no recycling mechanism is taken into account. Nevertheless the framework we set is very easy to manipulate: we can define and study other kinds of affinity-dependent selection mechanisms, and eventually include recycling mechanisms, which have been demonstrated to play an important role in AAM \cite{victora2010germinal}. We demonstrate that independently from the transition probability matrix defining the mutational mechanism and the affinity threshold chosen for positive selection, the optimal selection rate maximizing the number of output cells for the $t^{\rm th}$ generation is $1/t$ (Proposition \ref{cor:maxrs}), $t\in\mathbb N$.\\

From a mathematical point of view, we study a class of multi-types Galton-Watson (GW) processes (\emph{e.g.} \cite{harris2002theory,athreya2012branching}) in which, by considering dead and selected B-cells as two distinct types, we are able to formalize the evolution of a population submitted to an affinity-dependent selection mechanism. To our knowledge, the problem of affinity-dependent selection in GW processes has not been deeply investigated so far. \\

In Section \ref{sec0} we define the main model analyzed in this paper. We give as well some definitions that we will use in next sections. 
Section \ref{sec1} contains the main mathematical results. A convenient use of a multi-type GW process allows to study the evolution of both GC and output cells over time. 
We determine the optimal value of the selection rate which maximizes the expected number of selected B-cells at any given maturation cycle in Section \ref{sec:rsbest}.
We conclude Section \ref{sec1} with some numerical simulations. In Section \ref{sec:ext} we define two possible variants of the model described in previous sections, and provide some mathematical results and numerical simulations as well. 
This evidences how the mathematical tools used in Section \ref{sec1} easily apply to define other affinity-dependent selection models. Finally, in Section \ref{sec5} we discuss our modeling assumptions and give possible extensions and limitations of our mathematical model. In order to facilitate the reading of the paper, some technical mathematical demonstrations, as well as some classical results about Galton-Watson theory are reported in the Appendix for interested readers. \\

\section{Main definitions and modeling assumptions}\label{sec0}

This section provides the mathematical framework of this article. Let us suppose that given an antigen target cell $\cible$, all B-cell traits can be divided in exactly $N+1$ distinct affinity classes, named 0 to $N$.

\begin{definition}\label{def:affclass}
Let $\cible$ be the antigen target trait. Given a B-cell trait $\vec x$, we denote by $a_{\cible}(\vec x)$ the affinity class it belongs to with respect to  $\cible$, $a_{\cible}(\vec x)\,\in\,\{0,\dots,N\}$. The maximal affinity corresponds to the first class, 0, and the minimal one to $N$.
\end{definition}

\begin{definition}\label{def:aff}
Let $\vec x$ be a B-cell trait belonging to the affinity class $a_{\cible}(\vec x)$ with respect to  $\cible$. We say that its affinity with $\cible$ is given by:
\begin{equation*}
\textrm{aff}(\vec x,\cible)=N-a_{\cible}(\vec x)
\end{equation*}
\end{definition}

Of course, this is not the only possible choice of affinity. Typically affinity is represented as a Gaussian function \cite{wang2016stochastic,meyer2006analysis}, having as argument the distance between the B-cell trait and the antigen in the shape space of possible traits. In our model this distance corresponds to the index of the affinity class the B-cell belongs to (0 being the minimal distance, $N$ the maximal one). 
Nevertheless the choice of the affinity function does not affect our model.\\

During the GC reaction B-cells are submitted to random mutations. 
This implies switches from one affinity class to another  with a given probability. 
Setting these probability means defining a mutational rule on the state space $\{0,\dots,N\}$ of affinity classes indices (the formal mathematical definition will be given in Section \ref{sec3}).\\

\begin{figure}[ht!]
  \centering
\includegraphics[scale=0.4]{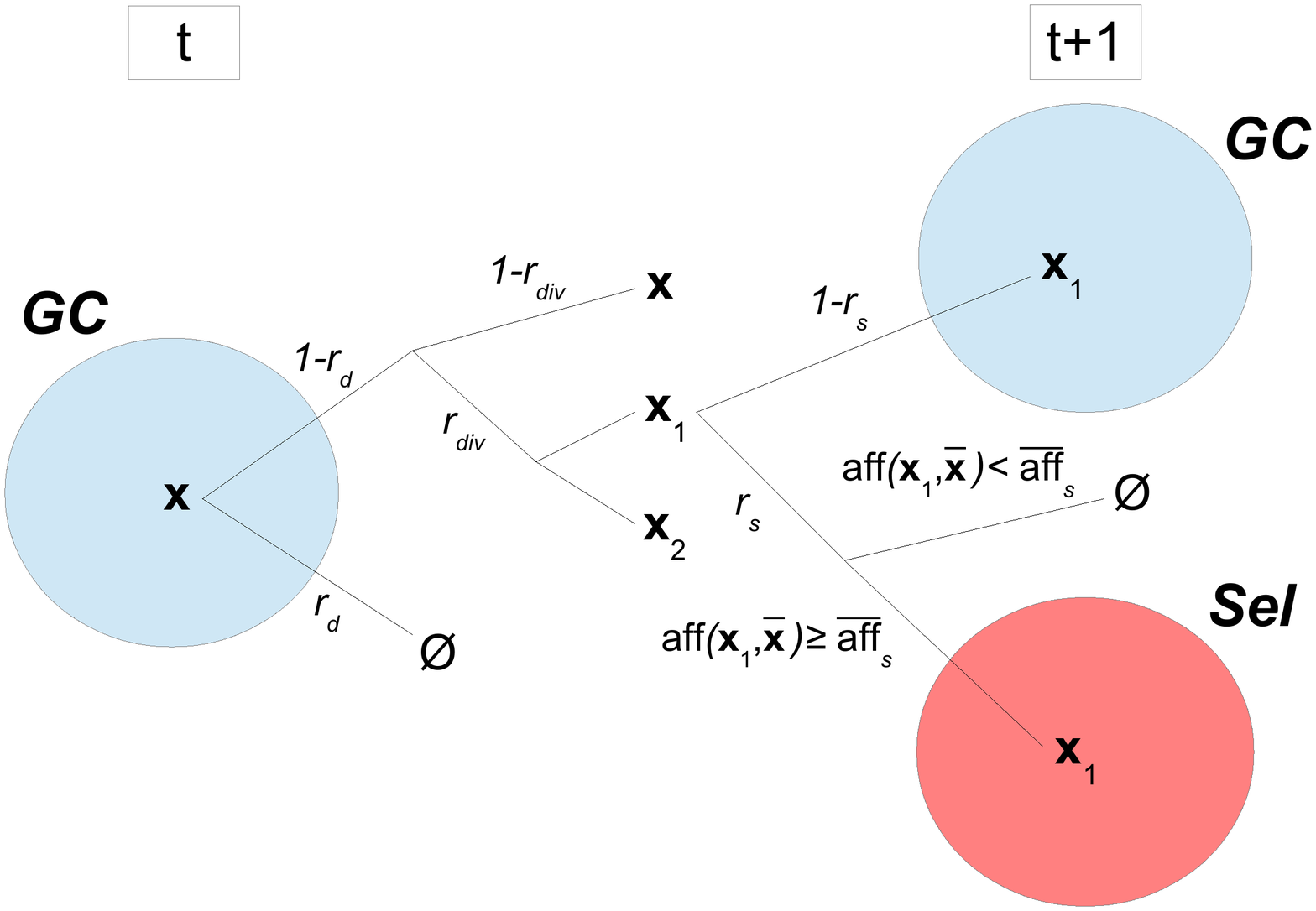}
  \caption{Schematic representation of model described by Definition \ref{def:basicmodel}. Here we denote by $\overline{\textrm{aff}}_s:=N-\overline a_s$, the fitness of each B-cell in the affinity class whose index is $\overline a_s$ (see Definitions \ref{def:affclass} and \ref{def:aff}).}
  \label{fig0}
\end{figure}

The main model we study in this paper is represented schematically in Figure \ref{fig0}. It is defined as follows: 
\begin{definition}\label{def:basicmodel}
The process starts with $z_0\geq1$ B-cells entering the GC, belonging to some affinity classes in $\{0,\dots,N\}$. In case they are all identical, we denote by $a_0$ the affinity class they belong to, with respect to  the antigen target cell $\cible$. At each time step, each GC B-cell can eventually undertake three distinct processes: division, mutation and selection. First of all, each GC B-cell can die with a given rate $r_d$. If not, each B-cell can divide with rate $r_{div}$: each daughter cell may have a mutated trait, according to the mutational rule allowed. Hence it eventually belongs to a different affinity class than its mother cell. 
Clearly, it  also happens that a B-cell stays in the GC without dying nor dividing. Finally, with rate $r_s$ each B-cell can be submitted to selection, which is made according to its affinity with $\cible$. A threshold $\overline{a}_s$ is fixed: if the B-cell belongs to an affinity class with index greater than $\overline{a}_s$, the B-cell dies. O\-ther\-wise, the B-cell exits the GC pool and reaches the selected pool. Therefore, for any GC B-cell and at any generation, we have:
\begin{equation*}
\left\{\begin{array}{llcl}
\textrm{Probability of cellular apoptosis:} & \Pro(\textrm{death}) & = & r_d \\
\textrm{Probability of cellular division:} & \Pro(\textrm{division}) & = & r_{div} \\
\textrm{Probability of selection challenge:} & \Pro(\textrm{selection}) & = & r_s \\
\end{array}\right.
\end{equation*}
\end{definition}

\captionsetup[figure]{labelfont=bf}

Once the GC reaction is fully established ($\sim$ day 7 after immunization), it is polarized into two compartments, named Dark Zone (DZ) and Light Zone (LZ) respectively. The DZ is characterized by densly packed dividing B-cells, while the LZ is less densely populated and contains FDCs and Tfh cells. The LZ is the preferential zone for selection \cite{de2015dynamics}. The transition of B-cells from the DZ to the LZ seems to be determined by a timed cellular program: over a 6 hours period about 50$\%$ of DZ B-cells transit to the LZ, where they compete for positive selection signaling \cite{bannard2013germinal,victora2014snapshot}. \\

Through the entire paper one should keep in mind the following main modeling assumptions:

\begin{modassumption}
In our simplified mathematical model we do not take into account any spatial factor and in a single time step a GC B-cell can eventually undergo both division (with mutation) and selection. Hence the time unit has to be chosen big enough to take into account both mechanisms.
\end{modassumption}

\begin{modassumption}
In this paper we are considering discrete-time models. The symbol $t$ always denote a discrete time step, hence it is an integral value. We will refer to $t$ as time, generation, or even maturation cycle to further stress the fact that in a single time interval $[t,t+1]$ each B-cell within the GC population is allowed to perform a complete cycle of division, mutation and selection.
\end{modassumption}

\begin{modassumption}
Throughout the entire paper, when we talk about death rate (respectively division rate or selection rate) we are referring to the probability that each cell has of dying (respectively dividing or being submitted to selection) in a single time step.\\
\end{modassumption}

\section{Results}\label{sec1}

In this Section we formalize mathematically the model introduced above. This enables the estimation of various qualitative and quantitative measures of the GC evolution and of the selected pool as well.
In Section \ref{sec2} we show that a simple GW process describes the evolution of the size of the GC and determine a condition for its extinction. 
In order to do this we do not need to know the mutational model. 
Nevertheless, if we want to understand deeply the whole reaction we need to consider a $(N+3)$-type GW process, which we introduce in Section \ref{sec3}. Therefore we determine explicitly other quantities, such as the average affinity in the GC and the selected pool, or the evolution of the size of the latter. We conclude this section by numerical simulations (Section \ref{sec4}).

\subsection{Evolution of the GC size}\label{sec2}

The aim of this section is to estimate the evolution of the GC size and its extinction probability. 
In order to do so we define a simple GW process, with respect to  parameters $r_{d}$, $r_{div}$ and $r_{s}$. Indeed, each B-cell submitted to selection exits the GC pool, independently from its affinity with $\cible$. Hence we apply some classical results about generating functions and GW processes (\cite{harris2002theory}, Chapter I), which we recall in Appendix \ref{app:classicGW}. 
Proposition \ref{prop:eta} gives explicitly the expected size of the GC at time $t$ and conditions for the extinction of the GC.

\begin{definition}\label{def:Zt}
Let $Z_t^{(z_0)}$, $t\geq0$ be the random variable (rv) describing the GC-population size at time $t$, starting from $z_0\geq1$ initial B-cells.
$(Z_t^{(z_0)})_{t\in\mathbb{N}}$ is a MC (as each cell behaves independently from the others and from previous generations) on $\{0,1,2,\dots\}$.
\end{definition}

If $z_0=1$ and there is no confusion, we denote $Z_t:=Z_t^{(1)}$. By Definition \ref{def:Zt}, $Z_1$ corresponds to the number of cells in the GC at the first generation, starting from a single seed cell. Thanks to Definition \ref{def:basicmodel} one can claim that $Z_1\,\in\,\{0,1,2\}$, with the following probabilities:

\begin{equation}\label{eq:proZt}
\left\{\begin{array}{l}
p_0:=\Pro(Z_1=0)=r_d+(1-r_d)r_s(1-r_{div}+r_{div}r_s) \\
p_1:=\Pro(Z_1=1)=(1-r_d)(1-r_s)(1-r_{div}+2r_{div}r_s)\\
p_2:=\Pro(Z_1=2)=r_{div}(1-r_d)(1-r_s)^2
\end{array}\right.
\end{equation}

As far as next generations are concerned, conditioning to $Z_t=k$, \emph{i.e.} at generation $t$ there are $k$ B-cells in the GC, $Z_{t+1}$ is distributed as the sum of $k$ independent copies of $Z_1$: $\Pro(Z_{t+1}=k'\,|\,Z_t=k)=\Pro\left(\sum_{i=1}^k Z_1=k'\right)$.\\

\begin{definition}
Let $X$ be an integer valued rv, $p_k:=\Pro(X=k)$ for all $k\geq0$. Its probability generating function (pgf) is given by:
\begin{equation*}
F_X(s)=\sum_{k=0}^{+\infty} p_k s^k
\end{equation*}
\end{definition}

The pgf for $Z_1$:

\begin{eqnarray}\label{eq:gen}
F(s) & = & p_0+p_1s+p_2s^2 \nonumber \\
& = & r_d+(1-r_d)r_s(1-r_{div}+r_{div}r_s) \nonumber \\
& + & (1-r_d)(1-r_s)(1-r_{div}+2r_{div}r_s)s +r_{div}(1-r_d)(1-r_s)^2s^2
\end{eqnarray}

By using classical results on Galton-Watson processes (see Appendix \ref{app:classicGW}), one can prove:

\begin{proposition}\label{prop:eta}
\begin{description}
\item[]
\item[(i)] The expected size of the GC at time $t$ and starting from $z_0$ initial B-cells is given by:
\begin{equation}\label{prop:expZt}
\E(Z_t^{(z_0)})=z_0\left((1-r_d)(1+r_{div})(1-r_s)\right)^{t}
\end{equation}
\item[(ii)] Denoted by $\eta_{z_0}$ the extinction probability of the GC population starting from $z_0$ initial B-cells, one has:
\begin{itemize}
\item if $r_s\geq1-\displaystyle\frac{1}{(1-r_d)(1+r_{div})}$, then $\eta_{z_0}=1$
\item otherwise $\eta_{z_0}=\eta^{z_0}<1$, $\eta$ being the smallest fixed point of \eqref{eq:gen}
\end{itemize}
\end{description}
\end{proposition}

In particular the process is subcritical or supercritical independently from $z_0$. 
In the supercritical case, increasing the number of B-cells at the beginning of the process makes the probability of extinction decrease. 
More precisely, in the case $\eta<1$, then $\eta_{z_0}\to0$ if $z_0\to\infty$, but we recall that GCs seem to be typically seeded by few B-cells, varying from ten to hundreds \cite{tas2016visualizing}.\\

\captionsetup[figure]{labelfont=bf}
\begin{figure}[ht!]
  \centering
\includegraphics[scale=0.3]{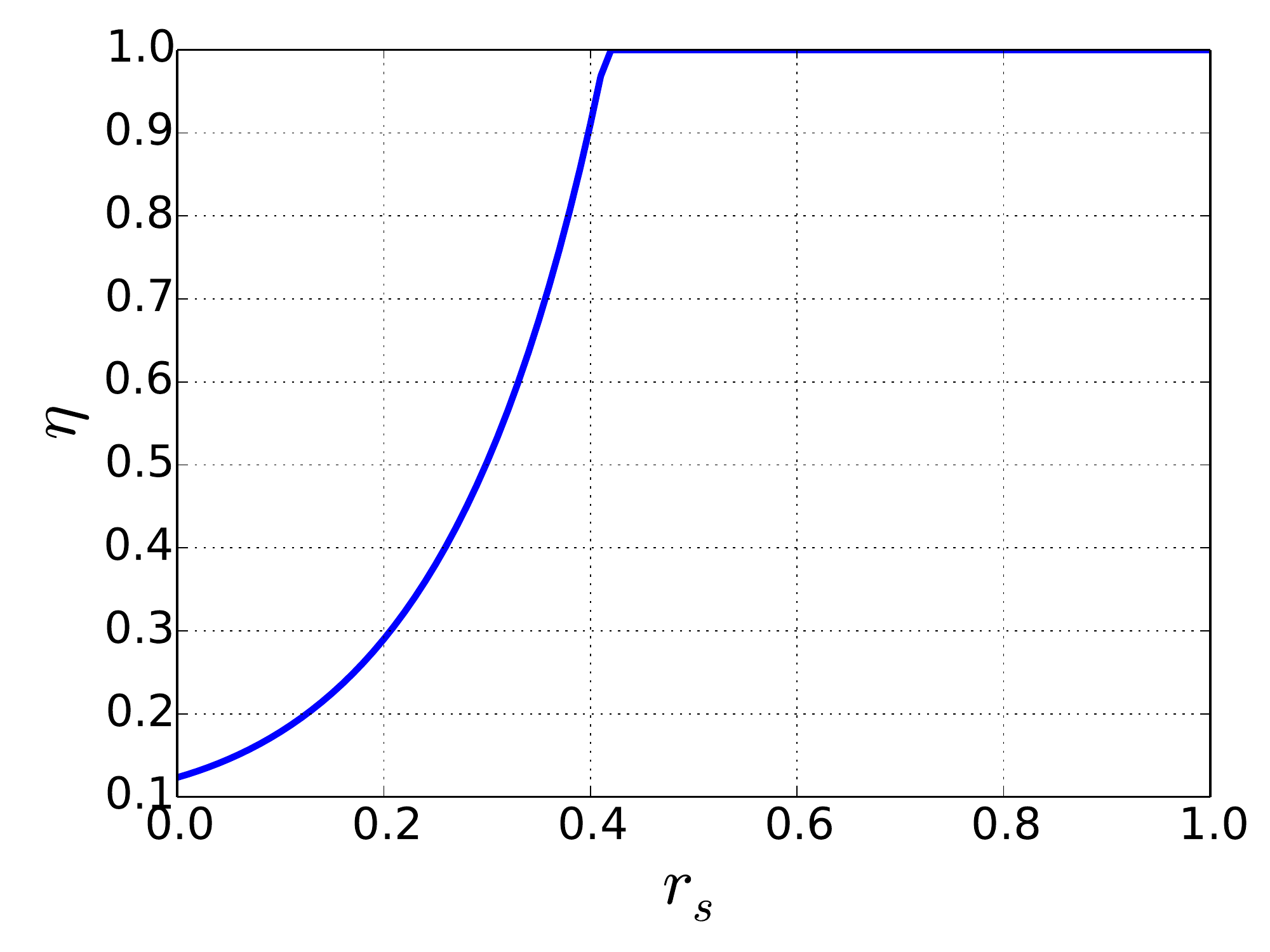}
  \caption{Numerical estimation of the extinction probability $\eta$ of the GC with respect to $r_s$ for $r_d=0.1$ and $r_{div}=0.9$.}
  \label{figextprob}
\end{figure}

This section shows that a classical use of a simple GW process enables to understand quantitatively the GC growth. Moreover, Proposition \ref{prop:eta} (ii) gives a condition on the main parameters for the extinction of the GC: if the selection pressure is too high, with probability 1 the GC size goes to 0, independently from the initial number of seed cells. Intuitively, a too high selection pressure prevents those B-cells with bad affinity to improve their fitness undergoing further rounds of mutation and division. Most B-cells will be rapidly submitted to selection, hence either exit the GC as output cells or die by apoptosis if they fail to receive positive selection signals \cite{maclennan2000b}. In Figure \ref{figextprob} we plot the extinction probability of a GC initiated from a single seed cell as a function of $r_s$ ($r_d$ and $r_{div}$ are fixed), in order to stress the presence of a threshold effect of the selection probability over the extinction probability. The extinction probability of the GC process can give us some further insights on factors which are potentially involved in determining the success or failure of a GC reaction. This simplified mathematical model suggests that if the selection pressure is too high compared to the division rate (\emph{c.f.} due to Tfh signals in the LZ), the GC will collapse with probability 1, preventing the generation of high affinity antibodies against the presented antigen, hence an efficient immune response.

\newcommand{\alert}[1]{{\color{red}{#1}}}

\subsection{Evolution of the size and fitness of GC and selected pools}\label{sec3}

The GW process defined in the previous Section only describes the size of the GC. Indeed, we are not able to say anything about the average fitness of GC clones, or the expected number of selected B-cells, or their average affinity. Hence, we need to consider a more complex model and take into account the parameter $\overline a_s$ and the transition probability matrix characterizing the mutational rule. 
Indeed, the mutational process is described as a Random Walk (RW) on the state space $\{0,\dots,N\}$ of affinity classes indices. The mutational rule reflects the edge set associated to the state-space $\{0,\dots,N\}$: this is given by a transition probability matrix.

\begin{definition}\label{def:QN}
Let $(\vec X_t)_{t\geq0}$ be a RW on the state-space of B-cell traits des\-cri\-bing a pure mutational process of a B-cell during the GC reaction. We denote by $\mathcal Q_N=(q_{ij})_{0\leq i,j\leq N}$ the transition probability matrix over $\{0,\dots,N\}$ which gives the probability of passing from an affinity class to another during the given mutational model. For all $0\leq i,j\leq N$:
\begin{equation*}
q_{ij}=\Pro(a_{\cible}(\vec X_{t+1})=j\,|\,a_{\cible}(\vec X_{t})=i)
\end{equation*}
\end{definition} 

We introduce a multi-type GW Process (see for instance \cite{athreya2012branching}, chapter V).

\begin{definition}\label{def:Ztmulti}
Let $\mathbf Z_t^{(\mathbf i)}=(Z_{t,0}^{(\mathbf i)},\dots,Z_{t,N+2}^{(\mathbf i)})$, $t\geq 0$ be a MC where for all $0\leq j\leq N$, $Z_{t,j}^{(\mathbf i)}$ describes the number of GC B-cells belonging to the $j^{\textrm{th}}$-affinity class with respect to  $\cible$, $Z_{t,N+1}^{(\mathbf i)}$ the number of selected B-cells and $Z_{t,N+2}^{(\mathbf i)}$ the number of dead B-cells at generation $t$, when the process is initiated in state $\mathbf i=(i_0,\dots,i_N,0,0)$. 
\end{definition}

Let $m_{ij}:=\E[Z_{1,j}^{(i)}]$ the expected number of offspring of type $j$ of a cell of type $i$ in one generation. We collect all $m_{ij}$ in a matrix, $\mathcal M=(m_{ij})_{0\leq i,j\leq N+2}$. We have:

\begin{equation}\label{eq:EZn}
\E[\mathbf Z_t^{(\mathbf i)}]=\mathbf i\mathcal M^t
\end{equation}

Supposing matrix $\mathcal Q_N$ given (Definition \ref{def:QN}), describing the probability to switch from one affinity class to another thanks to a single mutation event, one can explicitly derive the elements of $\mathcal M$.

\begin{proposition}\label{propM}
$\mathcal M$ is a $(N+3)\times(N+3)$ matrix defined as a block matrix:
\begin{equation*}
\mathcal M=\left(\begin{array}{cc}
\mathcal M_1 & \mathcal M_2 \\
\boldsymbol 0_{2\times (N+1)} & \mathcal I_2
\end{array}\right)
\end{equation*}
Where:
\begin{itemize}
\item $\boldsymbol 0_{2\times (N+1)}$ is a $2\times (N+1)$ matrix with all entries 0;
\item $\mathcal I_n$ is the identity matrix of size $n$;
\item $\mathcal M_1=2(1-r_d)r_{div}(1-r_s)\mathcal Q_N+(1-r_d)(1-r_{div})(1-r_s)\mathcal I_{N+1}$
\item $\mathcal M_2=(m_{2,ij})$ is a $(N+1)\times 2$ matrix where for all $i\,\in\,\{0,\dots,N\}$:
\begin{itemize}
\item if $i\leq\overline a_s$: \\
$m_{2,i1}=(1-r_d)(1-r_{div})r_s+2(1-r_d)r_{div}r_s\displaystyle\sum_{j=0}^{\overline a_s}q_{ij}$, \\
$m_{2,i2}=r_d+2(1-r_d)r_{div}r_s\displaystyle\sum_{j=\overline a_s+1}^{N}q_{ij}$
\item if $i>\overline a_s$: \\
$m_{2,i1}=2(1-r_d)r_{div}r_s\displaystyle\sum_{j=0}^{\overline a_s}q_{ij}$, \\
$m_{2,i2}=r_d+(1-r_d)(1-r_{div})r_s+2(1-r_d)r_{div}r_s\displaystyle\sum_{j=\overline a_s+1}^{N}q_{ij}$
\end{itemize}
\end{itemize}
\end{proposition}

The proof of Proposition \ref{propM} is available in Appendix \ref{app:proofPropM}. It is based on the computation of the probability generating function of $\mathbf Z_1$.

\begin{remark}
Independently from the given mutational model, the expected number of selected or dead B-cells that each GC B-cell can produce in a single time step is given by $\alpha:=r_d+(1-r_d)(1+r_{div})r_s$. All rows of $\mathcal M_2$ sum to $\alpha$ independently from the probability that each clone submitted to selection has of being positive selected, which we recall is 1 if it belongs to the $i^{\textrm{th}}$ affinity class, $i\leq\overline a_s$, zero otherwise.
\end{remark}

Of course in the multi-type context we recover again results from Section \ref{sec2}, such as the extinction probability of the GC (detailed in Appendix \ref{app:extinctionProbmulti}). \\

In order to determine the expected number of selected cells at a given time $t$, we need to introduce another multi-type GW process.

\begin{definition}
Let $\widetilde{\mathbf Z}_t^{(\mathbf i)}=(\widetilde{Z}_{t,0}^{(\mathbf i)},\dots,\widetilde{Z}_{t,N+2}^{(\mathbf i)})$, $t\geq 0$ be a MC where for all $0\leq j\leq N$, $\widetilde{Z}_{t,j}^{(\mathbf i)}$ describes the number of GC B-cells belonging to the $j^{\textrm{th}}$-affinity class with respect to  $\cible$, $\widetilde{Z}_{t,N+1}^{(\mathbf i)}$ the number of selected B-cells and $\widetilde{Z}_{t,N+2}^{(\mathbf i)}$ the number of dead B-cells at generation $t$, when the process is initiated in state $\mathbf i=(i_0,\dots,i_N,0,0)$ and before the selection mechanism is performed for the $t^{\textrm{th}}$-generation.
\end{definition} 

Proceeding as we did for $\mathbf Z_t^{(\mathbf i)}$, we can determine a matrix $\widetilde{\mathcal M}$ whose elements are $\widetilde{m}_{ij}:=\E[\widetilde{Z}_{1,j}^{(i)}]$ for all $i$, $j\;\in\;\{0,\dots,N+2\}$.

\begin{proposition}\label{propTildeM}
$\widetilde{\mathcal M}$ is a $(N+3)\times(N+3)$ matrix, which only depends on matrix $\mathcal Q_N$, $r_d$ and $r_{div}$ and can be defined as a block matrix as follows:
\begin{equation*}
\widetilde{\mathcal M}=\left(\begin{array}{cc}
\widetilde{\mathcal M}_1 & \widetilde{\mathcal M}_2 \\
\boldsymbol 0_{2\times (N+1)} & \mathcal I_2
\end{array}\right)
\end{equation*}
Where:
\begin{itemize}
\item $\widetilde{\mathcal M}_1=2(1-r_d)r_{div}\mathcal Q_N+(1-r_d)(1-r_{div})\mathcal I_{N+1}$
\item $\widetilde{\mathcal M}_2=\left(\boldsymbol{0}_{N+1},r_d\cdot\boldsymbol{1}_{N+1}\right)$, where $\boldsymbol{0}_{N+1}$ (resp. $\boldsymbol{1}_{N+1}$) is a $(N+1)$-column vector whose elements are all 0 (resp. 1).
\end{itemize}
\end{proposition}

One could prove that:
\begin{equation}\label{eq:ZtildeMGW}
\E \left[\widetilde{\mathbf Z}_t^{(\mathbf i)}\right]=\mathbf i\mathcal M^{t-1}\widetilde{\mathcal M}\
\end{equation}


\begin{proposition}\label{prop:expMM'}
Let $\mathbf i$ be the initial state, $|\mathbf i|$ its 1-norm ($|\mathbf i|:=\sum_{j=0}^{N+2}\mathbf i_j$).
\begin{itemize}
\item The expected size of the GC at time $t$:
\begin{equation}\label{eq:sizeGC}
\displaystyle\sum_{k=0}^N(\mathbf i\mathcal M^t)_k\left(=|\mathbf i|\left((1-r_d)(1+r_{div})(1-r_s)\right)^{t}\right)
\end{equation}
\item The average affinity in the GC at time $t$:
\begin{equation}\label{eq:affGC}
\displaystyle\frac{\displaystyle\sum_{k=0}^N(N-k)(\mathbf i\mathcal M^t)_k}{\displaystyle\sum_{k=0}^N(\mathbf i\mathcal M^t)_k}
\end{equation}
\item Let $S_t$, $t\geq1$ denotes the random variable describing the number of selected B-cells at time $t$. By hypothesis $S_0=0$.
$(S_t)_{t\in\mathbb{N}}$ is a MC on $\{0,1,2,\dots\}$.  The expected number of selected B-cells at time $t$, $t\geq1$:
\begin{equation}\label{eq:sizeSelt}
\E(S_t)=r_s\displaystyle\sum_{k=0}^{\overline a_s}\left(\mathbf i\mathcal M^{t-1}\widetilde{\mathcal M}\right)_k
\end{equation}
\item The expected number of selected B-cells produced until time $t$:
\begin{equation}\label{eq:sizeSelatt}
\E\left[\displaystyle\sum_{n=0}^tS_n\right]=\E\left[\left(\mathbf Z_t^{(\mathbf i)}\right)_{N+1}\right]=\left(\mathbf i\mathcal M^t\right)_{N+1}
\end{equation}
\item The average affinity of selected B-cells at time $t$, $t\geq1$:
\begin{equation}\label{eq:affSelt}
\displaystyle\frac{\displaystyle\sum_{k=0}^{\overline a_s}(N-k)\left(\mathbf i\mathcal M^{t-1}\widetilde{\mathcal M}\right)_k}{\displaystyle\sum_{k=0}^{\overline a_s}\left(\mathbf i\mathcal M^{t-1}\widetilde{\mathcal M}\right)_k}
\end{equation}
\item The average affinity of selected B-cells until time $t$:
\begin{equation}\label{eq:affSelatt}
\displaystyle\frac{r_s\displaystyle\sum_{s=1}^t\sum_{k=0}^{\overline a_s}(N-k)\left(\mathbf i\mathcal M^{s-1}\widetilde{\mathcal M}\right)_k}{\left(\mathbf i\mathcal M^t\right)_{N+1}}
\end{equation}
\end{itemize}
\end{proposition}

\proof
Equations \eqref{eq:sizeGC} and \eqref{eq:sizeSelatt} are a direct application of what stated in Equation \eqref{eq:EZn}. Indeed, Equation \eqref{eq:EZn} states that $\mathbf i\mathcal M^t$ contains the expectation of the number of all types cells at generation $t$ when the process is started in $\mathbf i$. Hence the expectation of the size of the GC at the $t^{\textrm{th}}$ generation is given by $\sum_{k=0}^N(\mathbf i\mathcal M^t)_k$, since the GC at generation $t$ contains all alive non-selected B-cells, irrespectively from their affinity. Similarly, the expected number of selected B-cells untill time $t$ \eqref{eq:sizeSelatt} corresponds to the expectation of the $(N+1)^{\textrm{th}}$-type cell, $\left(\mathbf i\mathcal M^t\right)_{N+1}$.\\

The proof of Equation \eqref{eq:sizeSelt} is based on Equation \eqref{eq:ZtildeMGW}, which allows to estimate the number of GC B-cells at generation $t$ which are susceptible of being challenged by selection. One can remark that the expected number of selected B-cells at time $t$ is obtained from the expected number of B-cells in GC at time $t$ (before the selection mechanism is performed) having fitness good enough to be positive selected. This is given by $\sum_{k=0}^{\overline a_s}\left(\mathbf i\mathcal M^{t-1}\widetilde{\mathcal M}\right)_k$, thanks to \eqref{eq:ZtildeMGW}. The result follows by multiplying this expectation by the probability that each of these B-cells is submitted to mutation, \emph{i.e.} $r_s$. Finally, results about the average affinity in both the GC and the selected pool (Equations \eqref{eq:affGC}, \eqref{eq:affSelt} and \eqref{eq:affSelatt}) are obtained from the pre\-vious ones (\emph{c.f.} \eqref{eq:sizeGC}, \eqref{eq:sizeSelt} and \eqref{eq:sizeSelatt}) by multiplying the number of individuals belonging to the same class by their fitness (Definition \ref{def:aff}), and dividing by the total number of individuals in the considered pool. The definition of affinity as a function of the affinity classes, determines Equations \eqref{eq:affGC}, \eqref{eq:affSelt} and \eqref{eq:affSelatt}.
\qed
\endproof

\begin{remark}
The expected size of the GC at time $t$ can be obtained applying a simple GW process (Section \ref{sec2}) and is given by  \eqref{prop:expZt}. 
It is possible to prove the equality in brackets in Equation \eqref{eq:sizeGC} starting from the  $(N+3)$-type GW process. The interested reader can address to Appendix \ref{app:proofeq:sizeGC} for the detailed proof.
\end{remark}

\subsection{Optimal value of $r_s$ maximizing the expected number of selected B-cells at time $t$}\label{sec:rsbest}

What is the behavior of the expected number of selected B-cells as a function of the model parameters ? In particular, is there an optimal value of the selection rate which maximizes this number ? In this section we show that, indeed, the answer is positive.\\

To do so we detail hereafter the computation  of $\E(S_t)$ (Equation \eqref{eq:sizeSelt}), given by Proposition \ref{prop:expMM'}. \\

Let us suppose, for the sake of simplicity, that $\mathcal Q_N$ is diagonalizable: 
\begin{equation}\label{eq:Qdiag}
\mathcal Q_N=R\Lambda_N L~,
 \end{equation}
 where $\Lambda_N=\diag(\lambda_0,\dots,\lambda_N)$, and $R=(r_{ij})$ (resp. $L=(l_{ij})$) is the transition matrix whose rows (resp. lines) contain the right (resp. left) eigenvectors of $\mathcal Q_N$, corresponding to $\lambda_0,\dots,\lambda_N$.  \\

\begin{proposition}\label{prop:expselatt}
Let us suppose that at $t=0$ there is a single B-cell entering the GC belonging to the $i^{\textrm{th}}$-affinity class with respect to  the target cell. Moreover, let us suppose that $\mathcal Q_N=R\Lambda_N L$. For all $t\in\mathbb N$, the expected number of selected B-cells at time $t$, is:
\begin{equation*}
\E(S_t)=r_s(1-r_s)^{t-1}(1-r_d)^t\displaystyle\sum_{\ell=0}^N(2\lambda_{\ell}r_{div}+1-r_{div})^t \sum_{k=0}^{\overline a_s}r_{i\ell}l_{\ell k}~,
\end{equation*}
\end{proposition}

The proof of Proposition \ref{prop:expselatt} is detailed in Appendix \ref{app:proofprop:expselatt}.\\

As an immediate consequence of Proposition \ref{prop:expselatt}, we can claim:

\begin{proposition}\label{cor:maxrs}
For all $t^{\ast}\in\mathbb N$ fixed, the value $r_s^{\ast}:=r_s(t^{\ast})$ which maximizes the expected number of selected B-cells at the ${t^{\ast}}^{\textrm{th}}$ maturation cycle is:
\begin{equation*}
r_s^{\ast}=\displaystyle\frac{1}{t^{\ast}}
\end{equation*}
\end{proposition}

\proof
Since $(1-r_d)^t\sum_{\ell=0}^N(2\lambda_{\ell}r_{div}+1-r_{div})^t\sum_{k=0}^{\overline a_s} r_{i\ell}l_{\ell k}$ is a non negative quantity independent from $r_s$, the value of $r_s$ which maximizes $\E(S_{t^{\ast}})$ is the one that maximizes $r_s(1-r_s)^{t^{\ast}-1}$. The result trivially follows.
\qed
\endproof

This result suggests that the selection rate in GCs is tightly related to the timing of the peak of a GC response. In particular, following this model, GCs which peak early (\emph{e.g.} for whom the maximal output cell production is reached in a few days) are possibly characterized by a higher selection pressure than GCs peaking later (the peak of a typical GC reaction has been measured to be close to day 12 post immunization \cite{wollenberg2011regulation}). Moreover, an high selection rate could also prevent a correct and efficient establishment of an immune response (\emph{c.f.} results about extinction probability - Proposition \ref{prop:eta}).\\

\begin{remark}
Under certain hypotheses about the mutational model and the GC evolution, one could justify the claim of Proposition \ref{cor:maxrs} by  heuristic arguments, without considering the $(N+3)$-type GW process. This leads to approximately estimate the expected number of selected B-cells at time $t$ (Appendix \ref{app:proofcor:maxrs}). 
\end{remark}

\subsection{Numerical simulations}\label{sec4}

We evaluate numerically results  of Proposition \ref{prop:expMM'}. 
The $(N+3)$-type GW process allows a deeper understanding of  the dynamics of both populations:
inside the GC and in the selected pool. Through numerical simulations 
we emphasize the dependence of the quantities defined in Proposition \ref{prop:expMM'}  on parameters involved in the model. \\

In previous works \cite{balelli2015branching,iba.2} we have modeled B-cells and antigens as $N$-length binary strings, hence their traits correspond to elements of $\{0,1\}^N$. In this context we have characterized affinity using the Hamming distance between B-cell and antigen representing strings. The idea of using a $N$-dimensional shape space to represent antibodies traits and their affinity with respect to a specific antigen has already been employed (\emph{e.g.} \cite{perelson1979theoretical,meyer2006analysis,kauffman1989nk}), and $N$ typically varies from 2 to 4. In the interests of simplification, we chose to set $N=2$. Moreover, from a biological viewpoint, this choice means that we classify the amino-acids composing B-cell receptors strings into 2 classes, which could represent amino-acids negatively and positively charged respectively. Charged and polar amino-acids are the most responsible in creating bonds which determine the antigen-antibody interaction \cite{KMMurPTraMWal}. \\

%
%

While performing numerical simulations (Sections \ref{sec4} and \ref{sec5:2}) we refer to the following transition probability matrix on $\{0,\dots,N\}$:

\begin{definition}\label{def:QPp}
For all $i$, $j\,\in\,\{0,\dots,N\}$:
\begin{equation*}
q_{ij}=\Pro(a_{\cible}(\mathbf X_{t+1})=j\,|\,a_{\cible}(\mathbf X_{t})=i)=\left\{\begin{array}{ll}
i/N & \textrm{if}\quad j=i-1 \\
(N-i)/N & \textrm{if}\quad j=i+1 \\
0 & \textrm{if}\quad |j-i|\neq 1 \\
\end{array}\right.
\end{equation*}
$\mathcal Q_N:=(q_{ij})_{0\leq i,j\leq N}$ is a tridiagonal matrix where the main diagonal consists of zeros.
\end{definition}

If we model B-cell traits as vertices of the state-space $\{0,1\}^N$, this corresponds to a model of simple point mutations (see \cite{balelli2015branching} for more details and variants of this basic mutational model on binary strings). 

\begin{example}\label{ex:M2pourP}
One can give explicitly the form of matrix $\mathcal M_2$ (Proposition \ref{propM}) corresponding to the mutational model defined in Definition \ref{def:QPp}:

\begin{equation*}
\mathcal M_2 =\;\left.\begin{array}{c}
\textcolor{gray50}{0} \\
\textcolor{gray50}{\vdots} \\
\textcolor{gray50}{\overline a_s-1} \\
\textcolor{gray50}{\overline a_s} \\
\textcolor{gray50}{\overline a_s+1} \\
\textcolor{gray50}{\overline a_s+2} \\
\textcolor{gray50}{\vdots} \\
\textcolor{gray50}{N}
\end{array}\right.\left(\begin{array}{cc}
 \alpha & r_d \\
\vdots & \vdots \\
 \alpha & r_d \\
 \alpha-\beta+\beta\frac{\overline a_s}{N} & r_d+\beta\frac{N-\overline a_s}{N} \\
 \beta\frac{\overline a_s+1}{N} & r_d+\alpha-\beta+\beta\frac{N-(\overline a_s+1)}{N} \\
0 & r_d+\alpha \\
\vdots & \vdots \\
0 & r_d+\alpha
\end{array}\right)~,
\end{equation*}
where:
\begin{itemize}
\item $\alpha:= (1-r_d)(1+r_{div})r_s$
\item $\beta:=2(1-r_d)r_{div}r_s$\\
\end{itemize}
\end{example}

\begin{remark}
Note that all mathematical results obtained in previous sections are independent from the mutation model defined in Definition \ref{def:QPp}. 
\end{remark}

We suppose that at the beginning of the process there is a single B-cell entering the GC belonging to the affinity class $a_0$.
Of course, the model we set allows to simulate any possible initial condition. Indeed, by fixing the initial vector $\vec i$, we can decide to start the reaction with more B-cells, in different affinity classes. When it is not stated otherwise, the employed parameter set for simulations is given in Table \ref{tab:par_sim}.\\

\begin{table}[h]
\caption{Parameter choice for simulations in Sections \ref{sec4} (unless stated otherwise).}
\label{tab:par_sim}       
\begin{center}
\begin{tabular}{c c c c c c }
\hline\noalign{\smallskip}
$\boldsymbol{N}$ & $\boldsymbol{r_s}$ & $\boldsymbol{r_d}$ & $\boldsymbol{r_{div}}$ & $\boldsymbol{a_0}$ & $\boldsymbol{\overline a_s}$  \\
\noalign{\smallskip}\hline\noalign{\smallskip} 
10 & 0.1 & 0.1 & 0.9 & 3 & 3 \\
\noalign{\smallskip}\hline
\end{tabular}
\end{center}
\end{table}

This parameter choice implies a small extinction probability (Proposition \ref{prop:eta}).

\subsubsection{Evolution of the GC population}

The evolution of the size of the GC can be studied by using the simple GW process defined in Section \ref{sec2}. Equation \eqref{prop:expZt}, in the case of a single initial B-cell, evidences that the expected number of B-cells within the GC for this model only depends on $r_d$, $r_{div}$ and $r_s$ and it is not driven by the initial affinity, nor by the threshold chosen for positive selection $\overline a_s$, nor by the  mutational rule. \\

Equation \eqref{prop:expZt} evidences that, independently from the transition pro\-ba\-bi\-li\-ty matrix defining the mutational mechanism, the GC size at time $t$ increases with $r_{div}$ and decreases for increasing $r_s$ and $r_d$. Moreover, the impact of these last two parameters is the same for the growth of the GC. One could expect this behavior since the effect of both the death and the selection on a B-cell is the exit from the GC.\\

In order to study the evolution of the average affinity within the GC, we need to refer to the $(N+3)$-type GW process defined in Section \ref{sec3}.

\captionsetup[figure]{labelfont=bf}
\begin{figure}[ht!]
\centering
\begin{tabular}{cc}
\subfloat[ ]{\label{fig4b}\includegraphics[width = 2.3in]{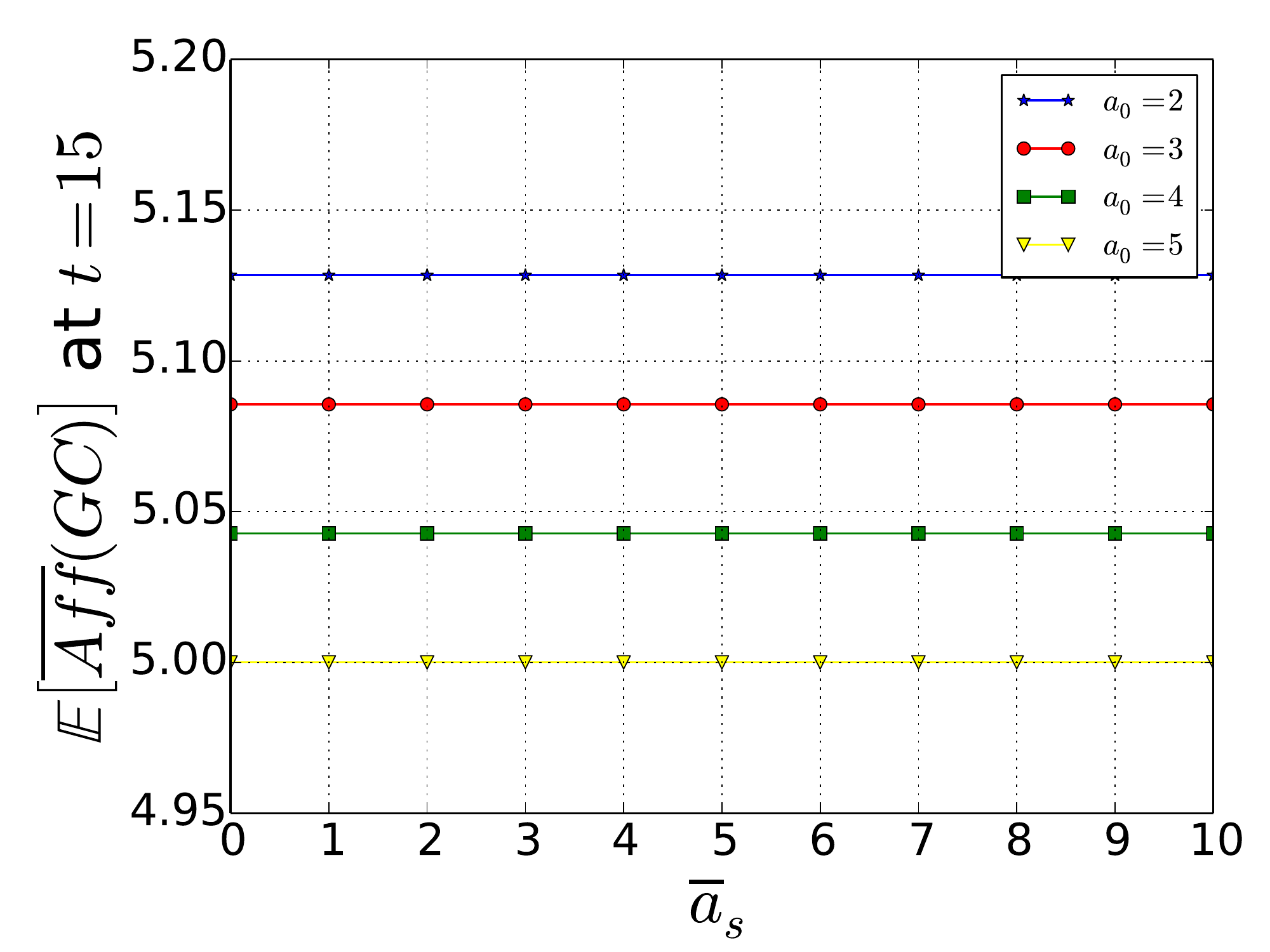}} & 
\subfloat[ ]{\label{fig4a}\includegraphics[width = 2.3in]{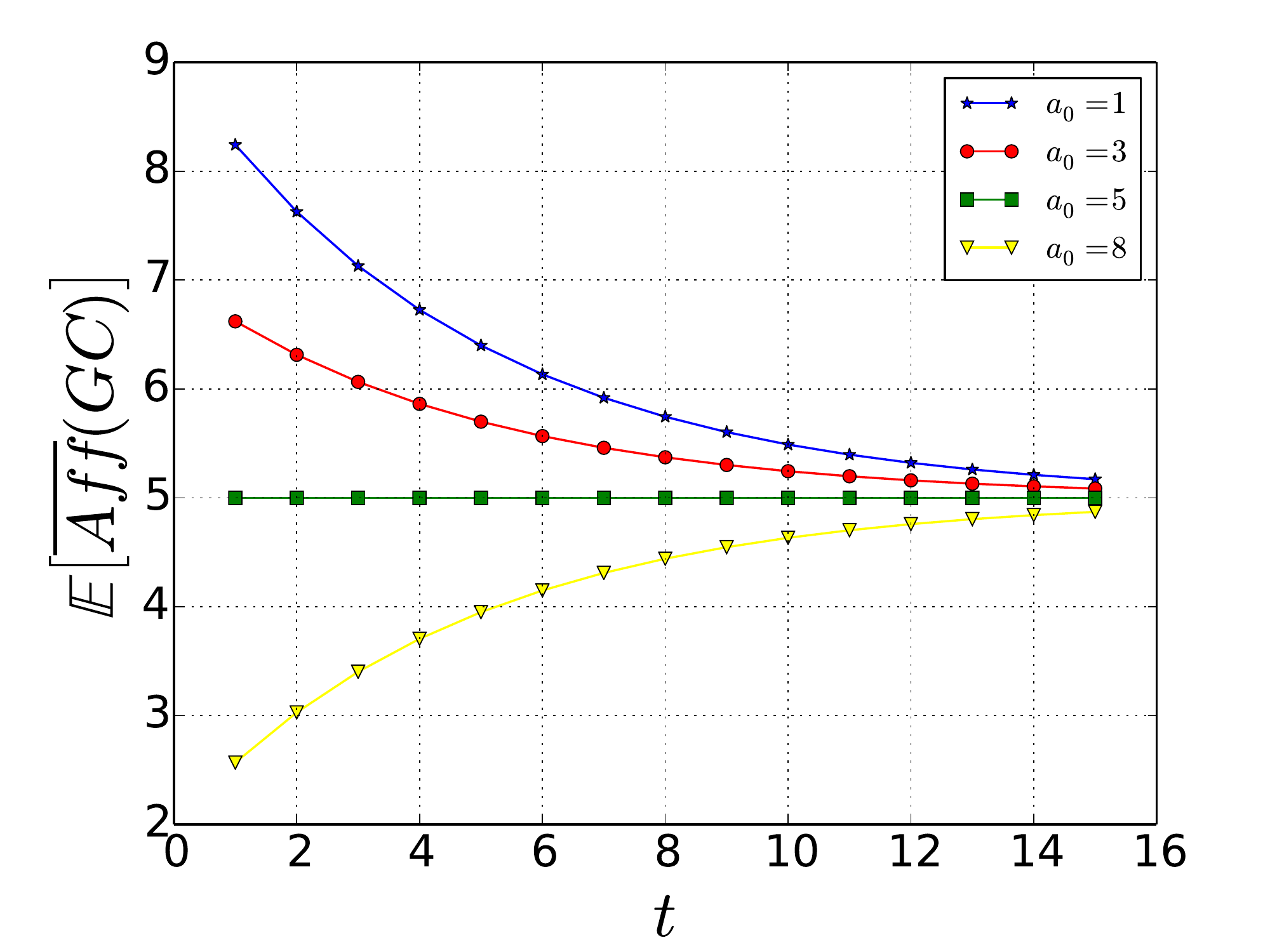}}\\
\end{tabular}
\caption{(a) Dependence of the expected average affinity in the GC on $\overline a_s$ at time $t=15$, for different values of $a_0$. The average affinity in the GC is constant with respect to $\overline a_s$. (b) The evolution during time of the expected average affinity in the GC for different values of $a_0$. The average affinity converges through $N/2$, due to the stationary distribution of $Q_N$, the binomial probability distribution.}\label{fig3_4}
\end{figure}

\begin{proposition}\label{prop:avGCM}
Let us suppose that $\mathcal Q_N=R\Lambda_N L$. The average affinity within the GC at time $t$, starting from a single B-cell belonging to the $i^{\textrm{th}}$-affinity class with respect to  $\cible$ is given by:
\begin{equation*}
N-\frac{\displaystyle\sum_{\ell=0}^N(2\lambda_{\ell}r_{div}+1-r_{div})^t \sum_{k=0}^Nk\cdot r_{i\ell}l_{\ell k}}{(1+r_{div})^t}~,
\end{equation*}
\end{proposition}

\proof
It follows directly from Equations \eqref{eq:affGC} and by considering the eigen\-de\-com\-po\-si\-tion of matrix $\mathcal Q$. One has to consider the expression of the $t^{\textrm{th}}$ power of matrix $\mathcal M$ (which can be obtained recursively, see Appendix \ref{app:proofprop:expselatt}): one can prove that the first $N+1$ components of the $i^{\textrm{th}}$-row of matrix $\mathcal M^t$ are the elements of the $i^{\textrm{th}}$-row of matrix $R D^t L$, where $D=2(1-r_d)r_{div}(1-r_s)\Lambda_N+(1-r_d)(1-r_{div})(1-r_s)\mathcal I_{N+1}$ is a diagonal matrix.
\qed
\endproof

It is obvious from Proposition \ref{prop:avGCM} that this quantity only depends on the initial affinity with the target trait, the transition probability matrix $\mathcal Q_N$ and the division rate $r_{div}$. The average affinity within the GC does not depend on  $\overline a_s$  (as one can clearly see in Figure \ref{fig3_4} (a)), nor by $r_s$ or $r_d$. One can intuitively understand this behavior:  independently from their fitness, all B-cells submitted to mutation exit the GC. 
Moreover, $r_s$ and $r_d$ impact the GC size, but not its average affinity, as selection and death affect all individuals of the GC independently from their fitness. \\

It can be interesting to observe the evolution of the expected average af\-fi\-ni\-ty within the GC during time. Numerical simulations of our model show that the expected average affinity in the GC converges through $N/2$, independently from the affinity of the first naive B-cell (Figure \ref{fig3_4} (b)). This depends on the mutational model we choose for these simulations.
Indeed, providing that the GC is in a situation of explosion, for $t$ big enough the distribution of GC clones within the affinity classes is governed by the stationary distribution of matrix $\mathcal Q_N$. Since for $\mathcal Q_N$ given by Definition \ref{def:QPp} one can prove that the stationary distribution over $\{0,\dots,N\}$ is the binomial probability distribution \cite{balelli2015branching}, the average affinity within the GC will quickly stabilizes at a value of $N/2$.

\subsubsection{Evolution of the selected pool}

The evolution of the number of selected B-cells during time necessarily depends on the evolution of the GC. In particular, let us suppose we are in the supercritical case, \emph{i.e.} the extinction probability of the GC is strictly smaller than 1. Than, with positive probability, the GC explodes and so does the selected pool. On the other hand, if the GC extinguishes, the number of selected B-cells will stabilize at a constant value, as once a B-cell is selected it can only stay unchanged in the selected pool.\\

\captionsetup[figure]{labelfont=bf}
\begin{figure}[hb!]
  \centering
  \begin{tabular}{cc}
  \subfloat[ ]{\label{fig9a}\includegraphics[scale=0.29]{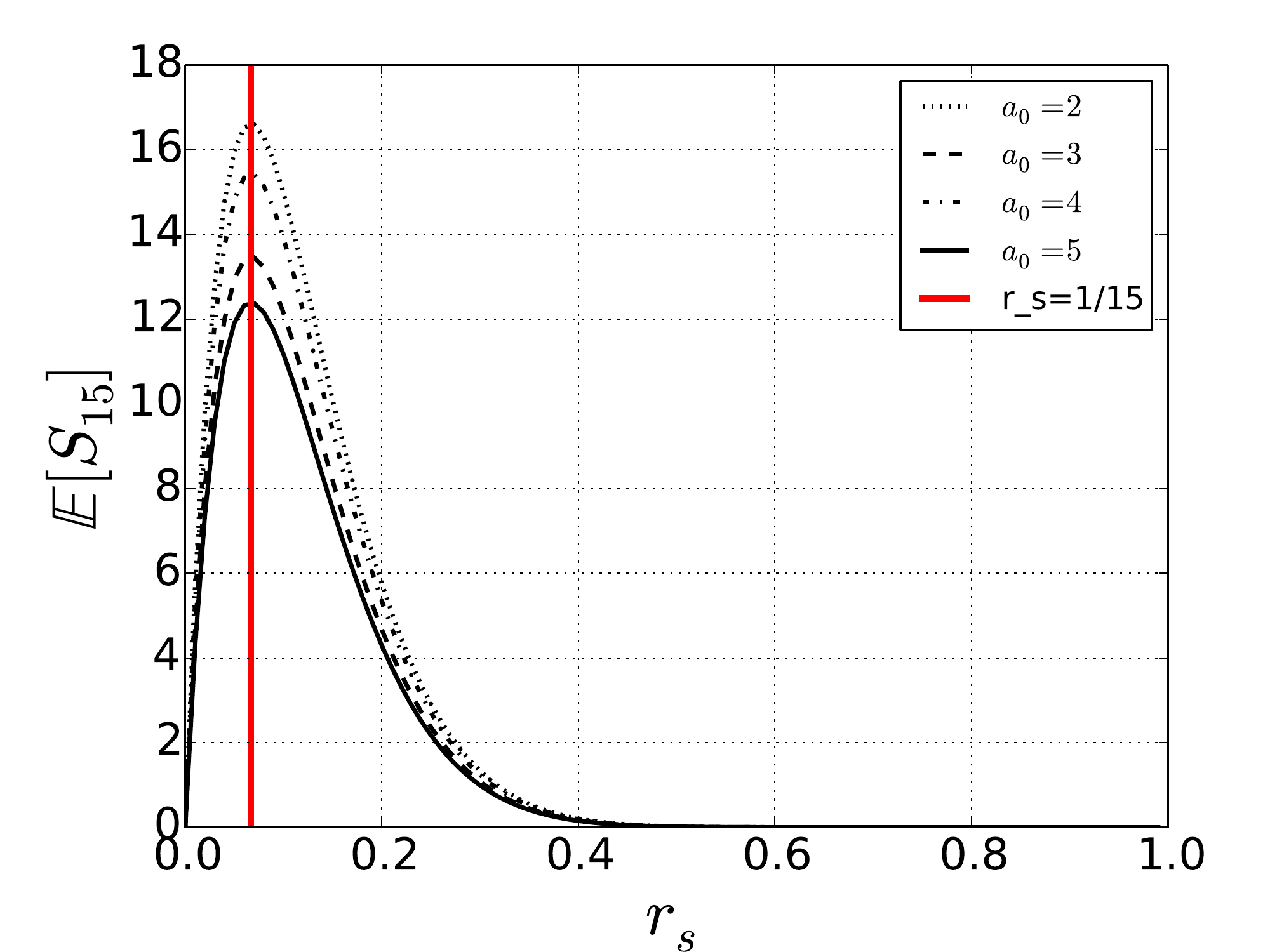}} 
  \subfloat[ ]{\label{fig9b}\includegraphics[scale=0.29]{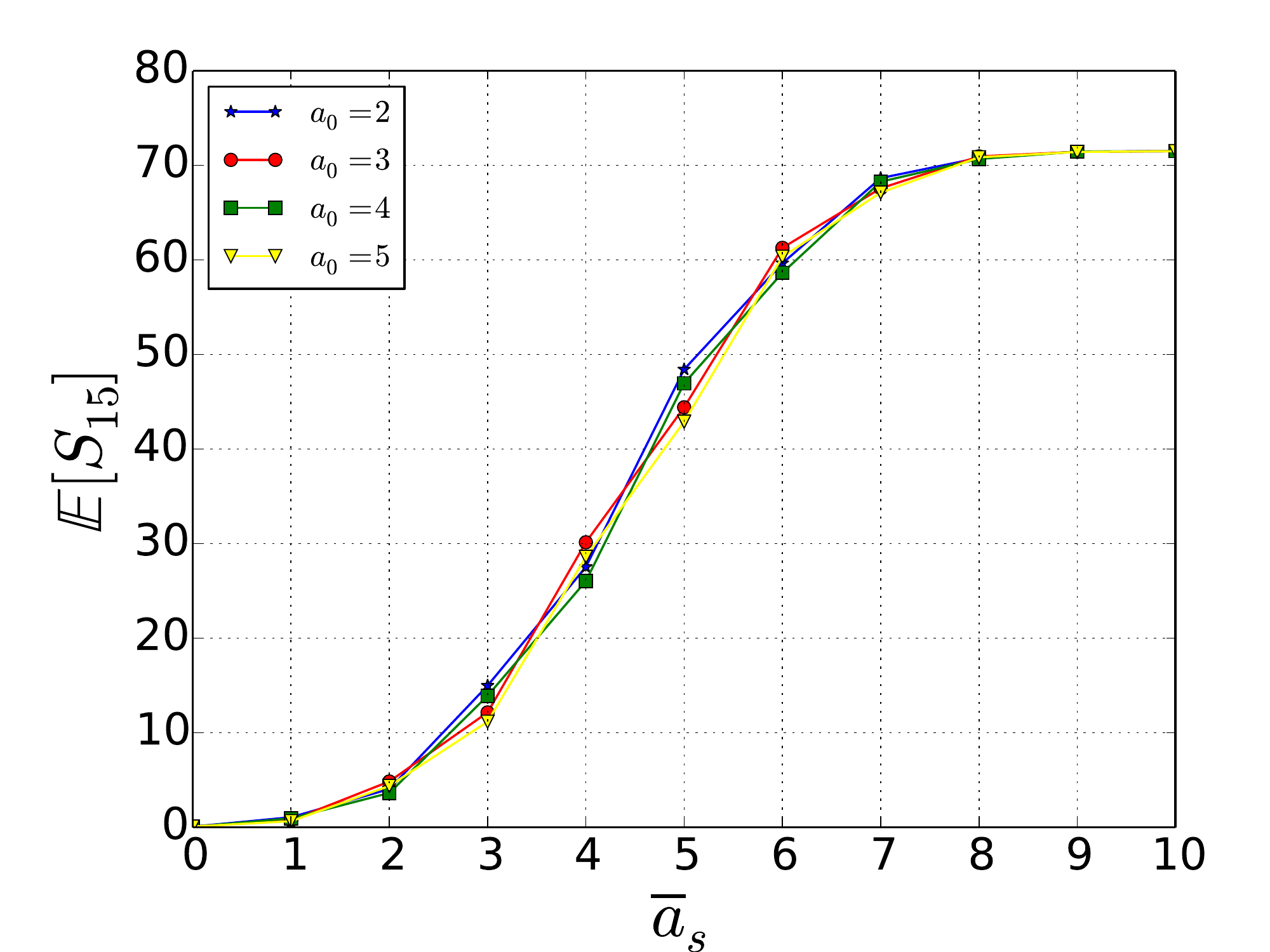}}\\
\end{tabular}
\subfloat[ ]{\label{fig10}\includegraphics[scale=0.29]{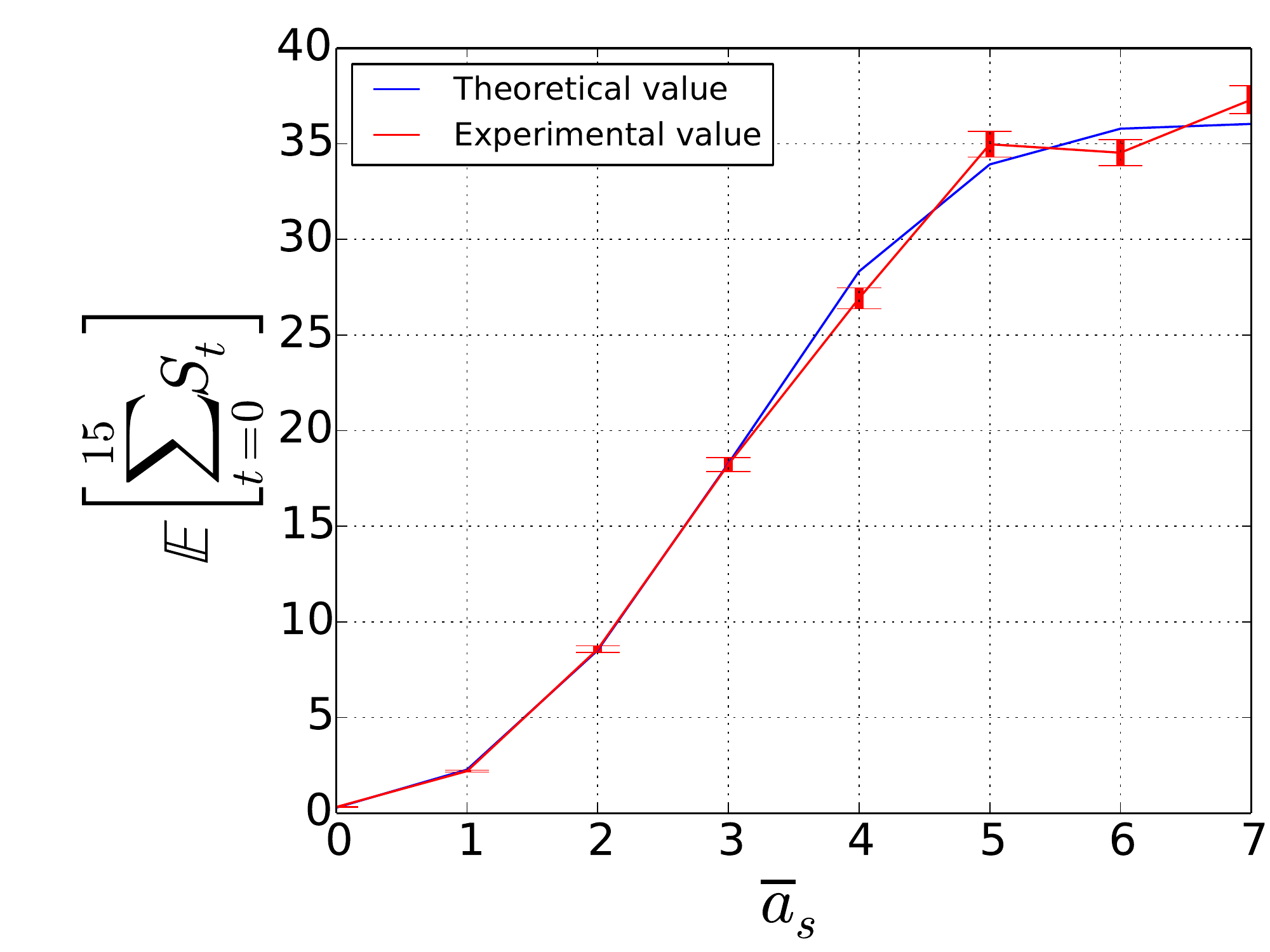}}
  \caption{(a-b) Expected number of selected B-cells for the time step $t=15$ for different values of $a_0$, depending on $r_s$ and $\overline a_s$ respectively. There exists an optimal value of $r_s$ maximizing the expected number of selected B-cells for a given generation. This value is independent from $a_0$ and is equal to $1/t$ as demonstrated in Proposition \ref{cor:maxrs}: the red vertical line in (a) corresponds to this value. (c) Comparison between the expected number of selected B-cells until time $t$ given by evaluation of the theoretical formula (Equation \eqref{eq:sizeSelatt}), and the empirical value obtained as the mean over 4000 simulations. Vertical bars denotes the corresponding estimated standard deviations. Here $N=7$ and $r_s=0.3$.}
  \label{fig9_10}
\end{figure}

As demonstrated in Section \ref{sec:rsbest}, there exists an optimal value of the parameter $r_s$ which maximizes the expected number of selected B-cells at time $t$. Figure \ref{fig9_10} (a) evidences this fact. Moreover, as expected, simulations show that the expected size of selected B-cells at a given time $t$ increases with the threshold $\overline a_s$ chosen for positive selection (Figure \ref{fig9_10} (b)). 
This  is a consequence  of Proposition \ref{prop:expselatt}: $\overline a_s$ determines the number of elements of the sum $\sum_{k=0}^{\overline a_s}r_{i\ell}l_{\ell k}$.\\

Figure \ref{fig9_10} (c) underlines the correspondence between  theoretical results given by Proposition \ref{prop:expMM'} and numerical values obtained by simulating the evolutionary process described by Definition \ref{def:basicmodel}. 
In particular Figure \ref{fig9_10} (c) shows the expected (resp. average) number of selected B-cells produced until time $t=15$ depending on the threshold chosen for positive selection, $\overline a_s$.

\section{Extensions of the model}\label{sec:ext}

Proceeding as in Section \ref{sec3}, we can define and study many different models of affinity-dependent selection. Here we propose a model in which we perform only positive selection and a model reflecting a Darwinian evolutionary system, in which the selection is only negative. For the latter, we will take into account only $N+2$ types instead of $N+3$: we do not have to consider a selected pool. Indeed the selected population remains in the GC. Here below we give the definitions of both models. In Section \ref{sec:resultsext} we formalize these problems mathematically, then in Section \ref{sec5:2} we show some numerical results.

\subsection{Definitions and results}\label{sec:resultsext}

Let us consider the process described in Definition \ref{def:basicmodel}. We change only the selection mechanism.

\captionsetup[figure]{labelfont=bf}
\begin{figure}[ht!]
  \centering
  \subfloat[Positive selection]{\label{fig11b}\includegraphics[scale=0.35]{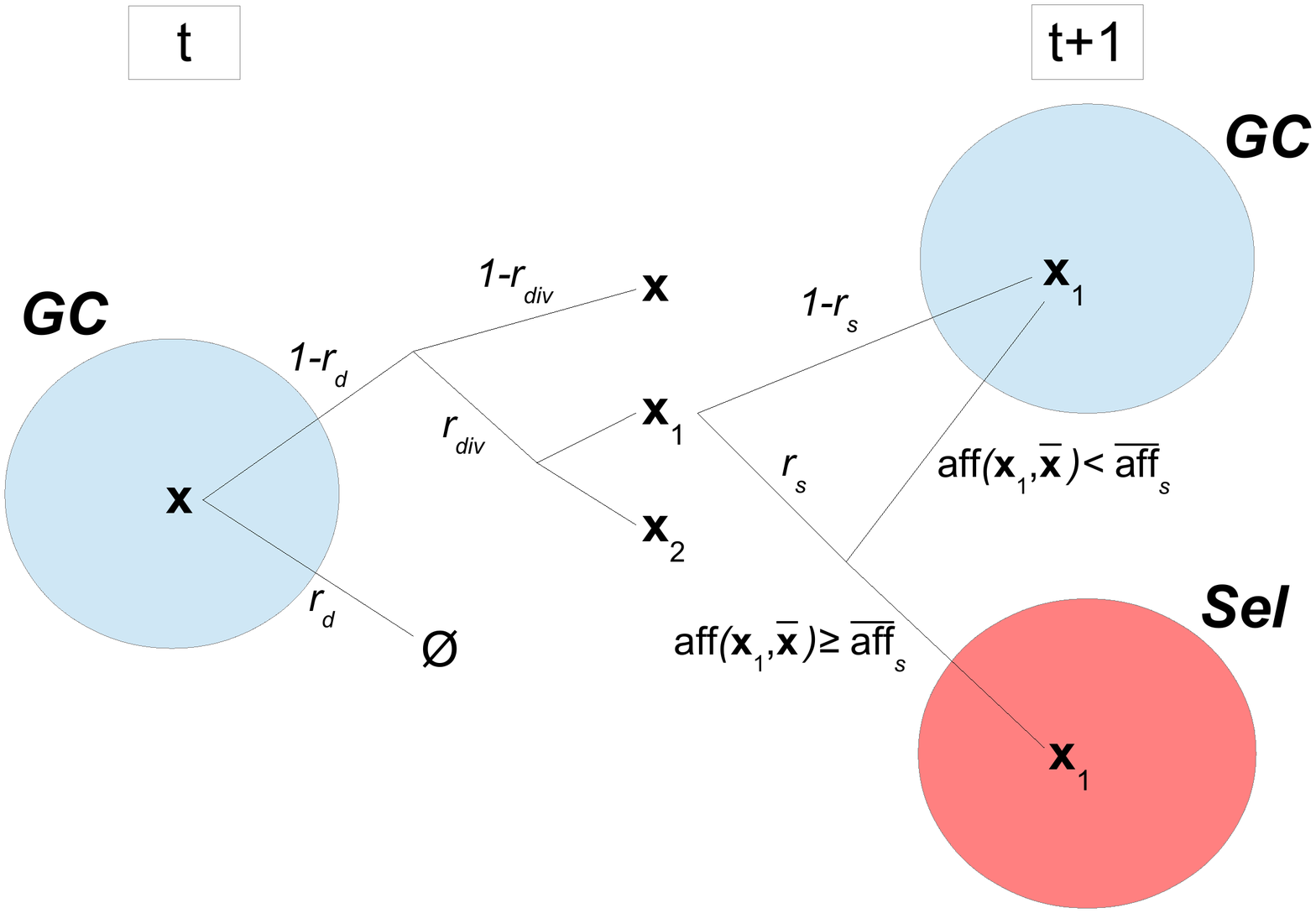}}
  \hspace{1pt}
  \subfloat[Negative selection]{\label{fig11a}\includegraphics[scale=0.35]{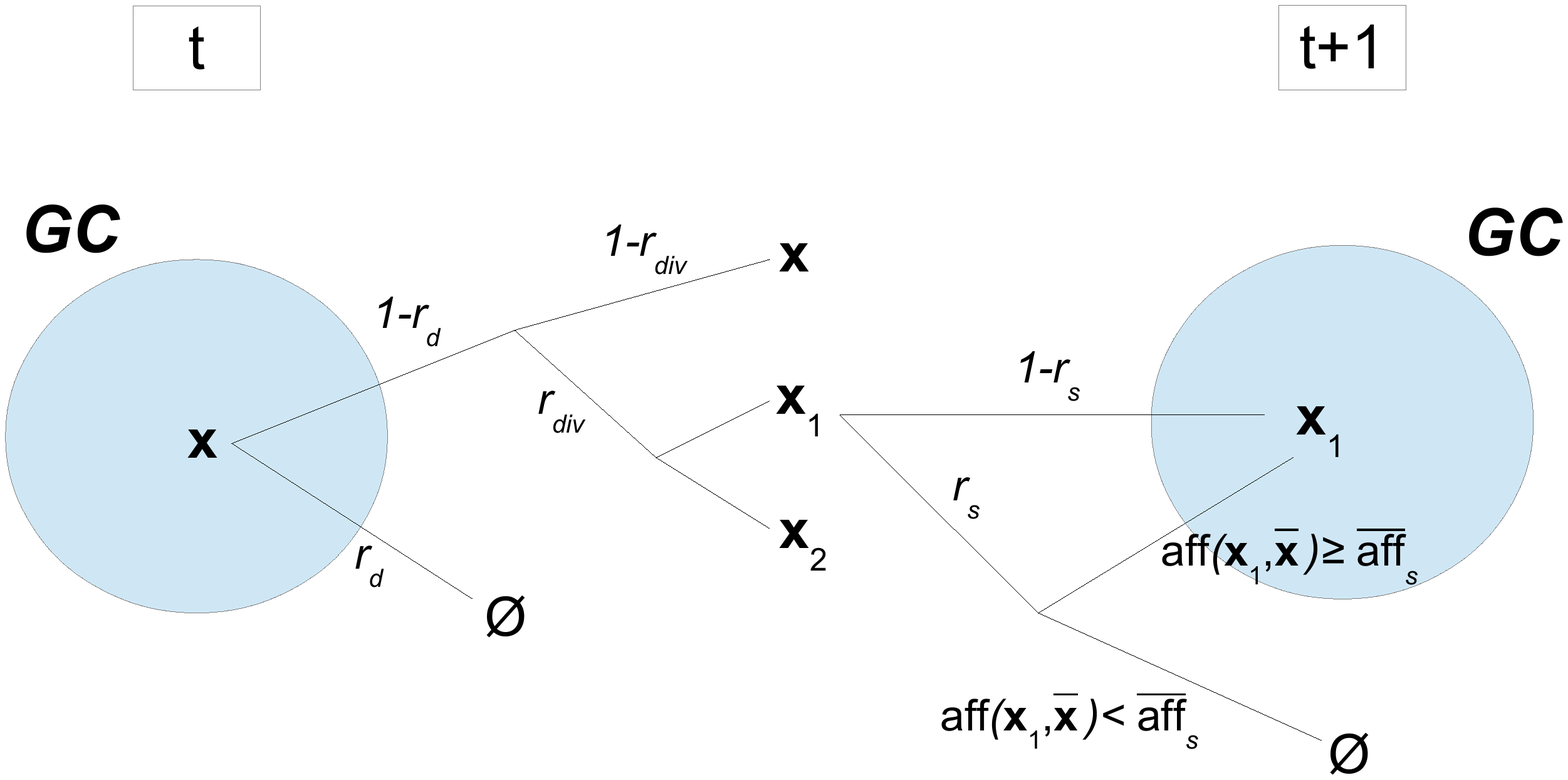}}
  \caption{Schematic representations of models described  (a) by Definitions \ref{def:modc} and (b) by Definitions \ref{def:moddarwin} of exclusively positive (resp. exclusively negative) selection.}
  \label{figpm}
\end{figure}

\begin{definition}[Positive selection]\label{def:modc}
If a B-cell submitted to selection belongs to an affinity class with index greater than $\overline{a}_s$, nothing happens. Otherwise, the B-cell exits the GC pool and reaches the selected pool.
\end{definition}

\begin{definition}[Negative selection]\label{def:moddarwin}
If a B-cell submitted to selection belongs to an affinity class with index greater than $\overline{a}_s$, it dies. Otherwise, nothing happens.
\end{definition}

In Figure \ref{figpm} we represent schematically both processes of positive selection and of negative selection. It is clear from Figure \ref{figpm} (b) that in the case of Definition \ref{def:moddarwin} we do not need to consider the selected pool anymore.

\paragraph{Positive selection}

\begin{definition}\label{def:Zt+}
Let ${\mathbf Z_t^+}^{(\mathbf i)}=({Z_{t,0}^+}^{(\mathbf i)},\dots,{Z_{t,N+2}^+}^{(\mathbf i)})$, $t\geq 0$ be a MC where for all $0\leq j\leq N$, ${Z_{t,j}^+}^{(\mathbf i)}$ describes the number of GC B-cells belonging to the $j^{\textrm{th}}$-affinity class with respect to  $\cible$, ${Z_{t,N+1}^+}^{(\mathbf i)}$ the number of selected B-cells and ${Z_{t,N+2}^+}^{(\mathbf i)}$ the number of dead B-cells at generation $t$, when the process is initiated in state $\mathbf i=(i_0,\dots,i_N,0,0)$, and following the evolutionary model described by Definition \ref{def:modc}.
\end{definition}

Let us denote by $\mathcal M^+=({m}_{ij}^+)_{0\leq i,j\leq N+2}$ the matrix containing the expected number of type-$j$ offsprings of a type-$i$ cell corresponding to the model defined by Definition \ref{def:modc}. We can explicitly write the value of all ${m}_{ij}^+$ depending on $r_d$, $r_{div}$, $r_s$, and the elements of matrix $\mathcal Q_N$.

\begin{proposition}\label{propM'}
$\mathcal M^+$ is a $(N+3)^2$ matrix, which we can define as a block matrix in the following way:
\begin{equation*}
\mathcal M^+=\left(\begin{array}{cc}
\mathcal {M}_1^+ & \mathcal {M}_2^+ \\
\boldsymbol 0_{2\times (N+1)} & \mathcal I_2
\end{array}\right)
\end{equation*}


Where:
\begin{itemize}
\item $\mathcal M^+_1=(m^+_{1,ij})$ is a $(N+1)^2$ matrix. For all $i\,\in\,\{0,\dots,N\}$:
\begin{itemize}
\item $\forall\,j\leq\overline a_s$: $m^+_{1,ij}=2(1-r_d)r_{div}(1-r_s)q_{ij}+(1-r_d)(1-r_{div})(1-r_s)\delta_{ij}$
\item $\forall\,j>\overline a_s$: $m^+_{1,ij}=2(1-r_d)r_{div}q_{ij}+(1-r_d)(1-r_{div})\delta_{ij}$
\end{itemize}
where $\delta_{ij}$ is the Kronecker delta.
\item $\mathcal M^+_2=(m^+_{2,ij})$ is a $(N+1)\times 2$ matrix where for all $i\,\in\,\{0,\dots,N\}$, $m^+_{2,i1}=m_{2,i1}$, and $m^+_{2,i2}=r_d$. We recall that $m_{2,i1}$ is the $i^{\textrm{th}}$-component of the first column of matrix $\mathcal M_2$, given in Proposition \ref{propM}.
\end{itemize}
\end{proposition}

\paragraph{Negative selection}

\begin{definition}\label{def:Ztmulti}
Let ${\mathbf Z_t^-}^{(\mathbf i)}=({Z_{t,0}^-}^{(\mathbf i)},\dots,{Z_{t,N+1}^-}^{(\mathbf i)})$, $t\geq 0$ be a MC where for all $0\leq j\leq N$, ${Z_{t,j}^-}^{(\mathbf i)}$ describes the number of GC B-cells belonging to the $j^{\textrm{th}}$-affinity class with respect to  $\cible$ and ${Z_{t,N+1}^-}^{(\mathbf i)}$ the number of dead B-cells at generation $t$, when the process is initiated in state $\mathbf i=(i_0,\dots,i_N,0)$, and following the evolutionary model described by Definition \ref{def:moddarwin}. 
\end{definition}

Let us denote by $\mathcal M^-=({m}_{ij}^-)_{0\leq i,j\leq N+1}$ the matrix containing the expected number of type-$j$ offsprings of a type-$i$ cell corresponding to the model defined by Definition \ref{def:Ztmulti}. 

\begin{proposition}\label{propM-}
$\mathcal M^-$ is a $(N+2)^2$ matrix, which we can define as a block matrix in the following way:
\begin{equation*}
\mathcal M^-=\left(\begin{array}{cc}
\mathcal {M}_1^- & \vec m_2^- \\
\boldsymbol 0_{N+1}' & 1
\end{array}\right)
\end{equation*}
Where:
\begin{itemize}
\item $\mathcal M^-_1=(m^-_{1,ij})$ is a $(N+1)^2$ matrix. For all $i\,\in\,\{0,\dots,N\}$:
\begin{itemize}
\item $\forall\,j\leq\overline a_s$: $m^-_{1,ij}=2(1-r_d)r_{div}q_{ij}+(1-r_d)(1-r_{div})\delta_{ij}$
\item $\forall\,j>\overline a_s$: $m^-_{1,ij}=2(1-r_d)r_{div}(1-r_s)q_{ij}+(1-r_d)(1-r_{div})(1-r_s)\delta_{ij}$
\end{itemize}
\item $\vec m^-_2$ is a $(N+1)$ column vector s.t. for all $i\,\in\,\{0,\dots,N\}$ $m^+_{i}=m_{2,i2}$, $m_{2,i2}$ being the $i^{\textrm{th}}$-component of the second column of matrix $\mathcal M_2$, given in Proposition \ref{propM}.
\item $\boldsymbol 0_{N+1}'$ is a $(N+1)$ row vector composing of zeros.
\end{itemize}
\end{proposition}

We do not prove Propositions \ref{propM'} and \ref{propM-}, since the proofs are the same as for  Proposition \ref{propM} (Appendix \ref{app:proofPropM}).\\

Results stated in Proposition \ref{prop:expMM'} hold true for these new models, by simply replacing matrix $\mathcal M$ with $\mathcal M^+$ (resp. $\mathcal M^-$). Of course, in the case of negative selection, as we do not consider the selected pool, we only refer to \eqref{eq:sizeGC} and \eqref{eq:affGC} quatifying the growth and average affinity of the GC. Matrix $\widetilde{\mathcal M}$ is the same for both models as only selection principles change.\\

Because of peculiar structures of  matrices $\mathcal M^+$ and $\mathcal M^-$, we are not able to compute explicitly their spectra. 
Henceforth we can not give an explicit formula for the extinction probability or evaluate the optimal values of the selection rate $r_s$ as we did in Sections \ref{sec3} and \ref{sec:rsbest}. \\

Nevertheless, by using standard arguments for positive matrices, 
  the grea\-test eigenvalue of both matrices $\mathcal M_1^+$ and $\mathcal M_1^-$ can be bounded, and hence give  sufficient conditions for extinction. Indeed, form classical results about multi-type GW processes, the value of the greatest eigenvalue allows to discriminate between subcritical case (\emph{i.e.} extinction probability equal to 1) and supercritical case (\emph{i.e.} extinction probability strictly smaller than 1) \cite{athreya2012branching}.

\begin{proposition}\label{prop:estextpm}
Let $\vec q^+$ (resp. $\vec q^-$) be the extinction probability of the GC for the model corresponding to matrix $\mathcal M_1^+$ (resp. $\mathcal M_1^-$). 
\begin{itemize}
\item If $r_{div}\leq\displaystyle\frac{r_d}{1-r_d}$, then $\vec q^+=\vec q^-=\boldsymbol 1$. \\
\item If $r_s<1-\displaystyle\frac{1}{(1-r_d)(1+r_{div})}$, then $\vec q^+<\boldsymbol 1$ and $\vec q^-<\boldsymbol 1$.
\end{itemize}
\end{proposition}

\proof
Since both matrices $\mathcal M_1^+$ and $\mathcal M_1^-$ are strictly positive matrices, the Perron Frobenius 
Theorem insures that the spectral radius is also the greatest eigenvalue. Then the following classical result holds \cite{minc1988nonnegative}:

\begin{theorem}\label{thm:boundspecrad}
Let $A=(a_{ij})$  be a square nonnegative matrix with spectral radius $\rho(A)$ and let $r_i(A)$ denote the sum of the elements along the $i^{\textrm{th}}$-row of $A$. Then:
\begin{equation*}
\min_i r_i(A)\leq\rho(A)\leq\max_i r_i(A)
\end{equation*}
\end{theorem}

Simple calculations provide:
$$
\begin{aligned}
 \min_i r_i(\mathcal M_1^+)&=(1-r_d)(1+r_{div})-r_s(1-r_d)\left(2r_{div}\min_i\sum_{j=0}^{\overline a_s}q_{ij}+1-r_{div}\right) \\
\max_i r_i(\mathcal M_1^+)&=(1-r_d)(1+r_{div})-2r_sr_{div}(1-r_d)\max_i\sum_{j=0}^{\overline a_s}q_{ij} \\
\min_i r_i(\mathcal M_1^-)&=(1-r_d)(1+r_{div})-r_s(1-r_d)\left(2r_{div}\min_i\sum_{j=\overline a_s+1}^{N}q_{ij}+1-r_{div}\right)\\
\max_i r_i(\mathcal M_1^-)&=(1-r_d)(1+r_{div})-2r_sr_{div}(1-r_d)\max_i\sum_{j=\overline a_s+1}^{N}q_{ij}\\
\end{aligned}
$$

The result follows by observing that for all $i\,\in\,\{0,\dots,N\}$, $0\leq\sum_{j=0}^{\overline a_s}q_{ij}$,
$\sum_{j=\overline a_s+1}^{N}q_{ij}\leq1$, and applying Theorem \ref{thm:extmultit}.
\qed
\endproof

 \captionsetup[figure]{labelfont=bf}
\begin{figure}[ht!]
  \centering
 \includegraphics[scale=0.3]{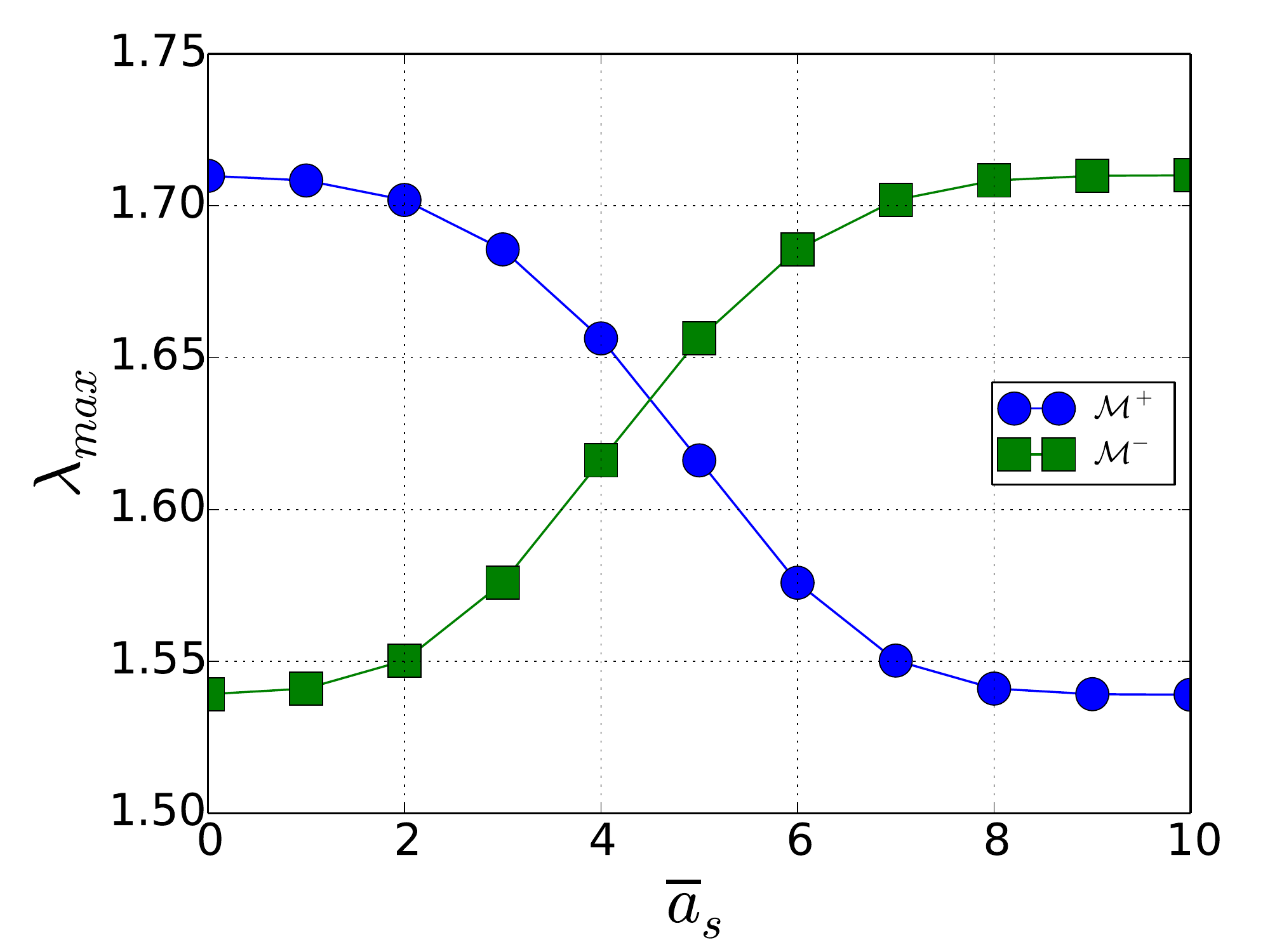}
  \caption{Dependence of greater eigenvalues of matrices $\mathcal M^+$ (blue circles) and $\mathcal M^-$ (green squares) respectively on $\overline a_s$ for $N=10$, $r_{div}=0.9$, $r_d=r_s=0.1$. Hence $(1-r_d)(1+r_{div})(1-r_s)=1.539$ and $(1-r_d)(1+r_{div}) =1.71$.}
  \label{figlammax}
\end{figure}

\begin{remark}
One can intuitively obtain the second claim of Proposition \ref{prop:estextpm}, as this condition over the parameters implies that the probability of extinction of the GC for the model underlined by matrix $\mathcal M_1$ of positive and negative selection is strictly smaller than 1 (Proposition \ref{prop:eta}). Indeed keeping the same parameters for all models, the size of the GC for the model of positive and ne\-ga\-tive selection is smaller than the size of GCs corresponding to both models of only positive and only negative selection. Consequently if the GC corresponding to $\mathcal M$ has a positive probability of explosion, it will be necessarily the same for $\mathcal M^+$ and $\mathcal M^-$.
\end{remark}

\begin{remark}\label{remmaxlam}
The values of both $\rho(\mathcal M_1^+)$ and $\rho(\mathcal M_1^-)$ depend on $\overline a_s$, varying from a minimum of $(1-r_d)(1+r_{div})(1-r_s)$ and a maximum of $(1-r_d)(1+r_{div})$.
Figure \ref{figlammax} evidences the dependence on $\overline a_s$ of the spectral radius of $\mathcal M_1^+$ and $\mathcal M_1^-$, 
using matrix $\mathcal Q_N$ given by Definition \ref{def:QPp} as transition probability matrix.
\end{remark}

Remark \ref{remmaxlam} and Figure \ref{figlammax} evidences that, conversely to the previous case of positive and negative selection, in both cases of exclusively positive (resp. exclusively negative) selection the parameter $\overline a_s$ plays an important role in the GC dynamics, affecting its extinction probability. In particular, keeping unchanged all other parameters, if $\overline a_s\to N$ (resp. $\overline a_s\to 0$), then $\rho(\mathcal M_1^+)$ (resp. $\rho(\mathcal M_1^-)$) $\to(1-r_d)(1+r_{div})(1-r_s)$, which implies $\vec q^+$ (resp. $\vec q^-$) $\to\boldsymbol 1$.\\

 \captionsetup[figure]{labelfont=bf}
\begin{figure}[ht!]
  \centering
 \includegraphics[scale=0.3]{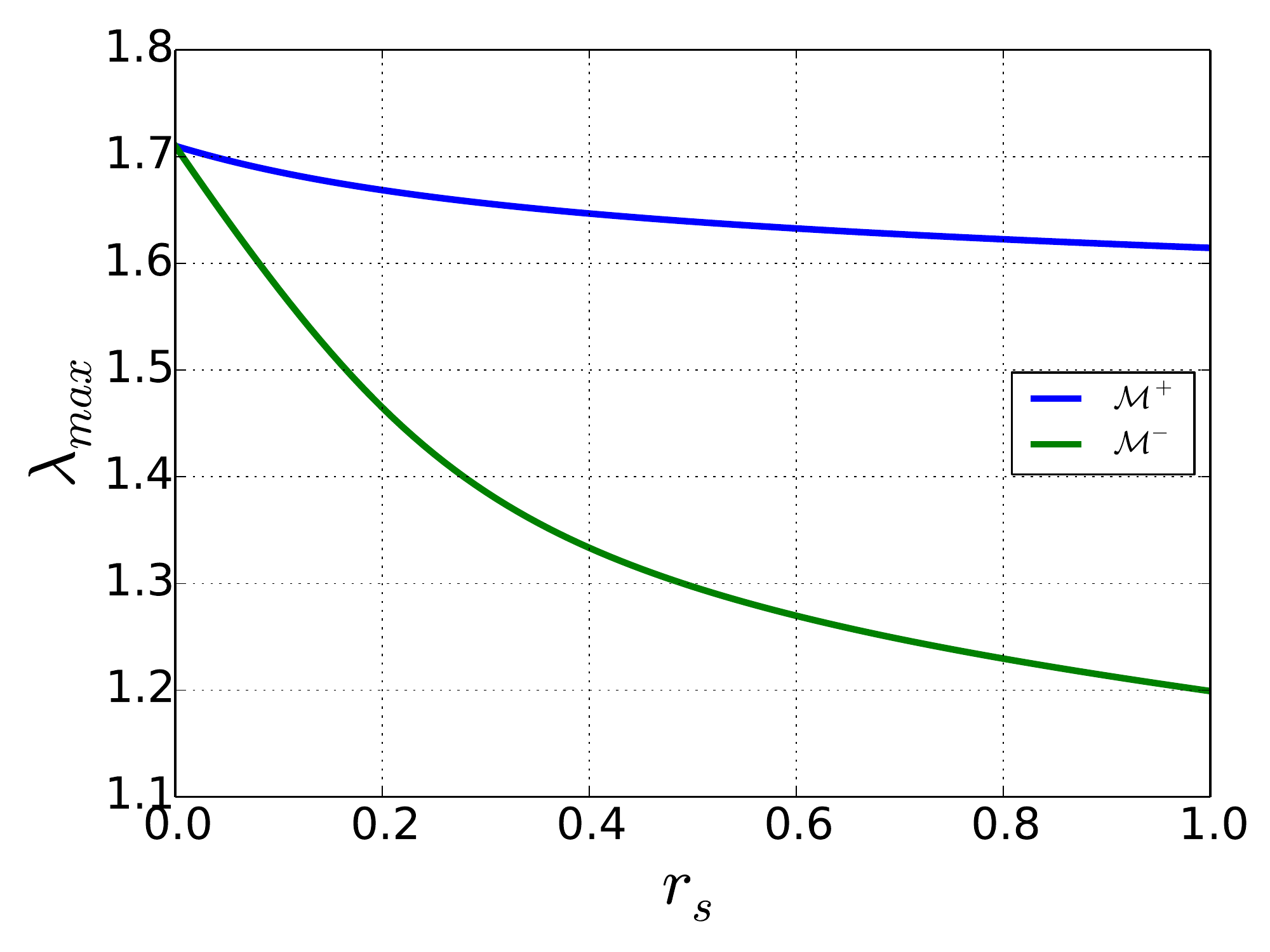}
  \caption{Dependence of greater eigenvalues of matrices $\mathcal M^+$ (blue) and $\mathcal M^-$ (green) respectively on $r_s$ for $N=10$, $r_{div}=0.9$, $r_d=0.1$, $\overline a_s=3$.}
  \label{figlammaxrs}
\end{figure}

In Figure \ref{figlammaxrs} we plot the dependence of greater eigenvalues of both matrices $\mathcal M^+$ and $\mathcal M^-$ with respect to $r_s$. We fix $r_d=0.1$ and $r_{div}=0.9$ as for Figure \ref{figextprob}. One can note that with this parameter set and if $\overline a_s$ is chosen not ``too small'' nor ``too large'' with respect to $N$, then the greater eigenvalue for both matrices is always greater than 1 independently from $r_s$, \emph{i.e.} the extinction probability is always strictly smaller than 1.\\

\subsection{Numerical simulations}\label{sec5:2}

\captionsetup[figure]{labelfont=bf}
\begin{figure}[ht!]
\centering
\begin{tabular}{cc}
 \subfloat[ ]{\label{fig113}\includegraphics[width = 2.3in]{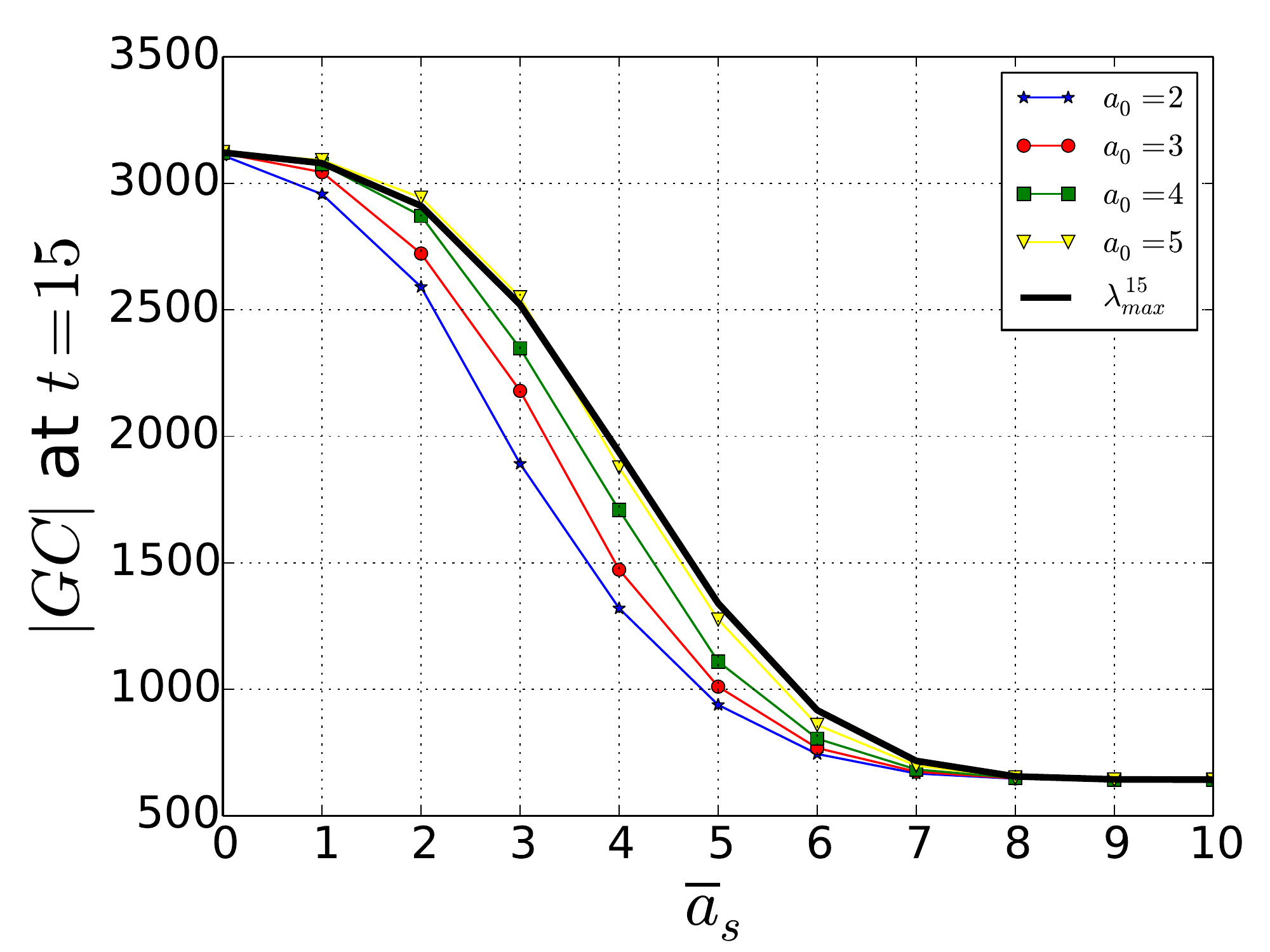}} &
\subfloat[ ]{\label{fig112}\includegraphics[width = 2.3in]{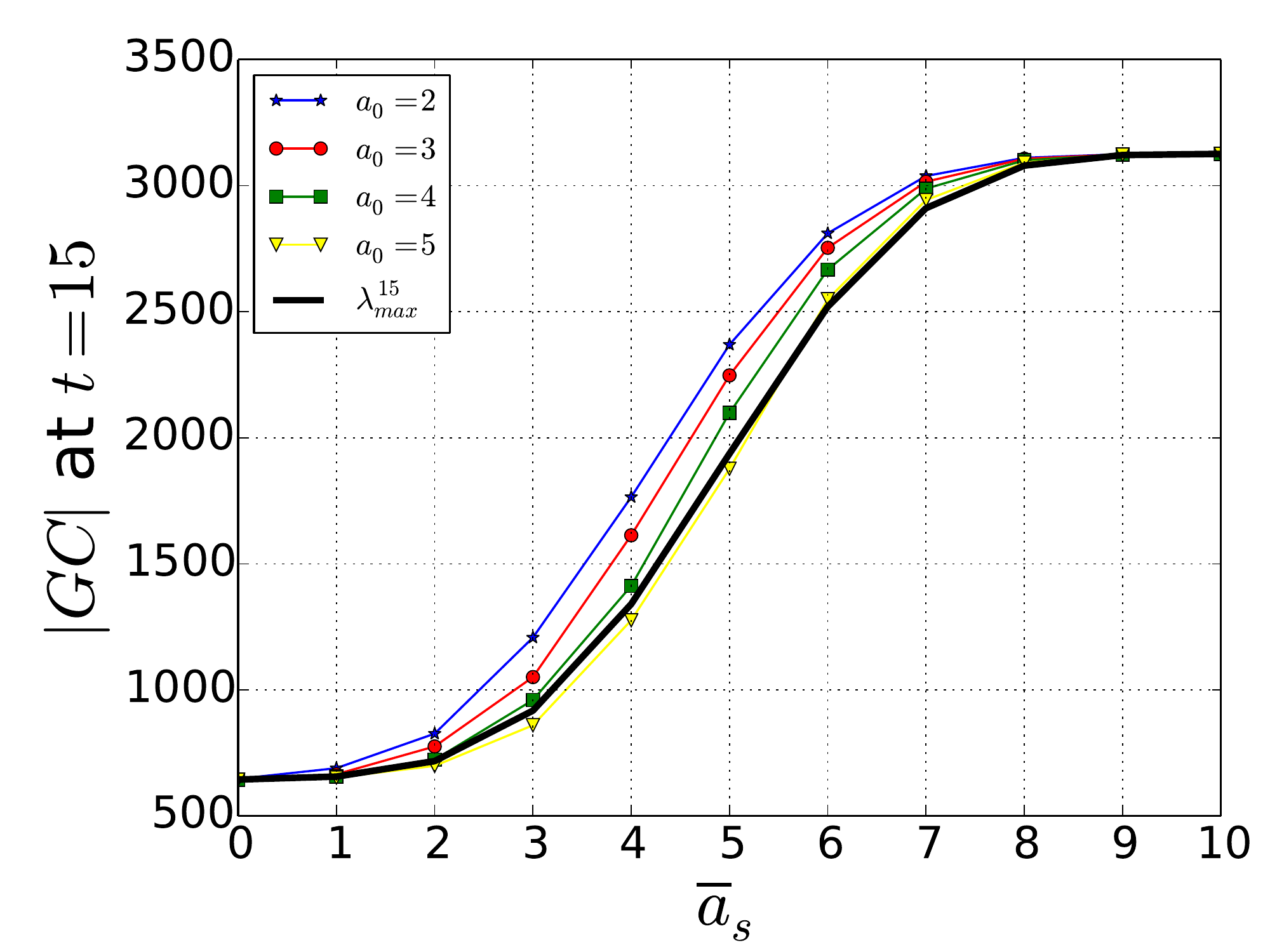}} \\
\subfloat[ ]{\label{fig114a}\includegraphics[width = 2.3in]{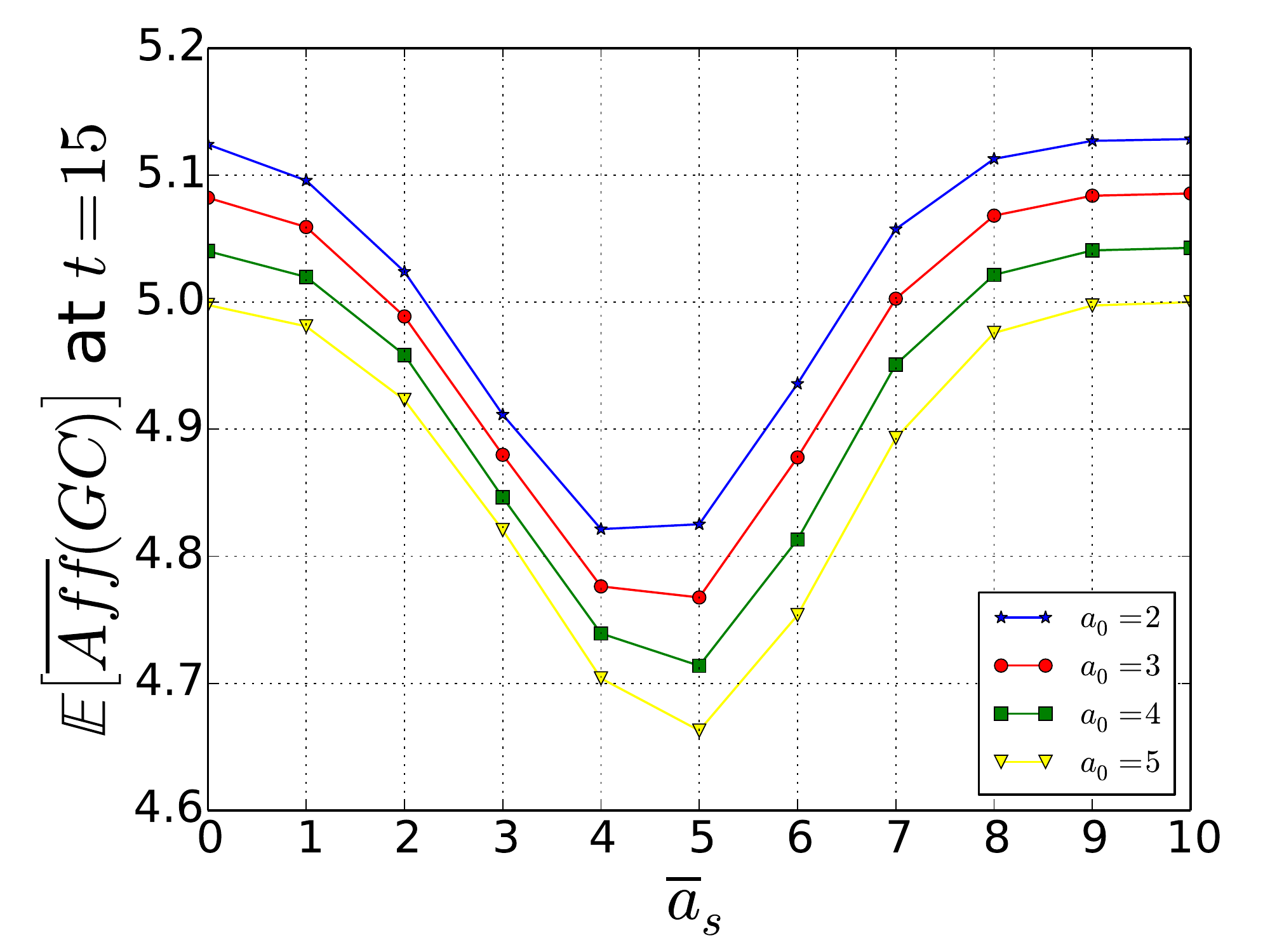}} & 
\subfloat[ ]{\label{fig114b}\includegraphics[width = 2.3in]{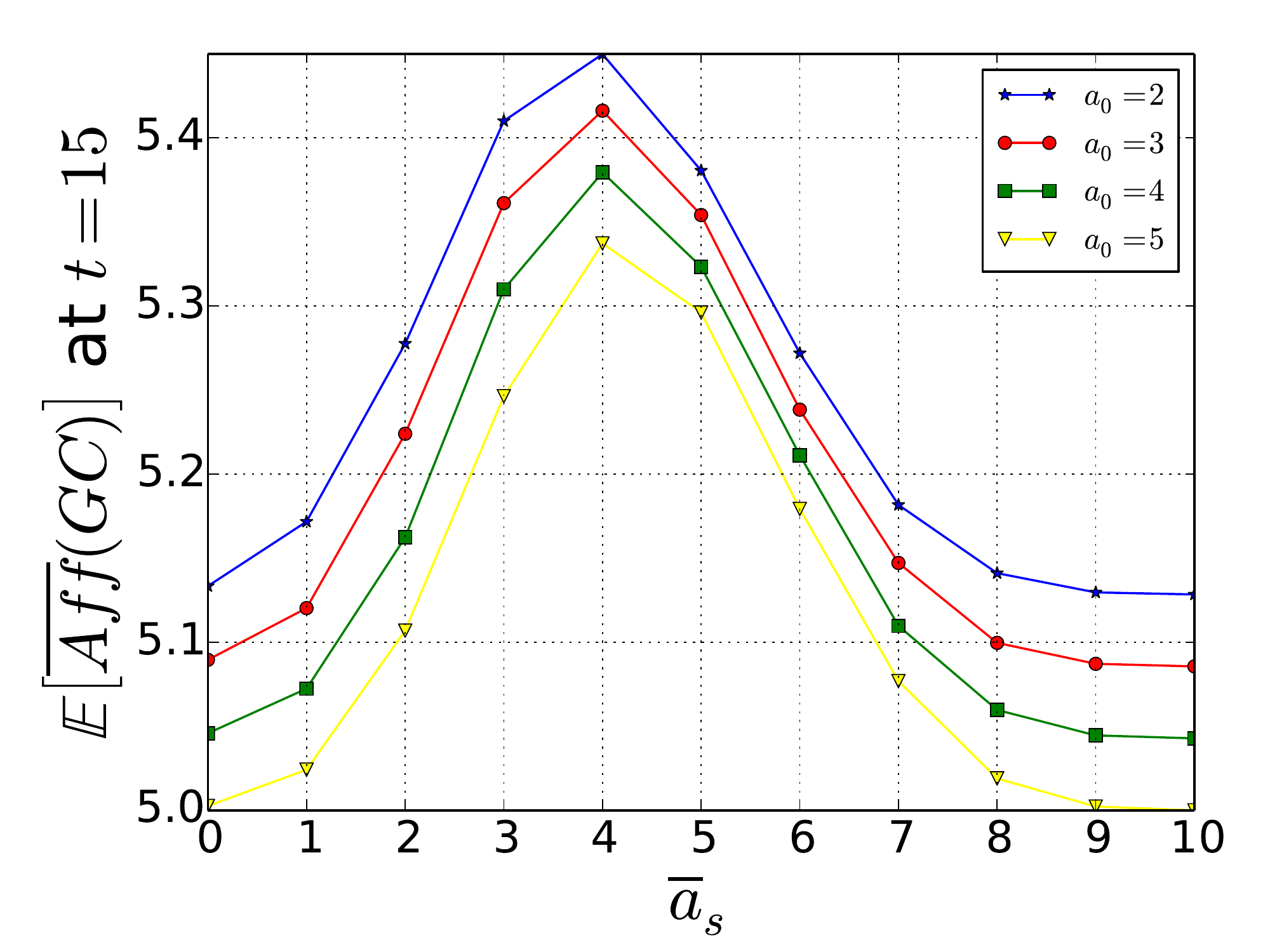}} \\
\end{tabular}
\caption{(a,b) Dependence of the expected size of the GC after 15 time steps on $\overline a_s$ for different values of $a_0$. The thick black line corresponds in both figures to the value of the greater eigenvalue of matrices $\mathcal M_1^+$ and $\mathcal M_1^-$ respectively, raised to the power of $t=15$ (see Figure \ref{figlammax}). Note that thanks to Proposition \ref{prop:eta} we know that for this parameter choice the expected size of the GC for the model of positive and negative selection corresponds to $((1-r_d)(1+r_{div})(1-r_s))^{15}$, which is equivalently $\lambda_{max}^{15}$ for $\overline a_s=10$ in Figure \ref{fig5_11} (a) or $\lambda_{max}^{15}$ for $\overline a_s=0$ in Figure \ref{fig5_11} (b). (c,d) Dependence of the expected average affinity in the GC after $t=15$ time steps on $\overline a_s$ for different values of $a_0$. The left column of Figure \ref{fig5_11} refers to the model of positive selection, while the right column to the model of negative selection.}\label{fig5_11}
\end{figure}

The evolution of GCs corresponding to matrices $\mathcal M^+$ and $\mathcal M^-$ respectively are complementary. 
Moreover, in both cases, keeping all parameters fixed one  expects a faster expansion if compared to the model of positive and negative selection, 
since the selection acts only positively (resp. negatively) on good (resp. bad) clones.
In particular, the model of negative selection 
corresponds to the case of 100\% of recycling, meaning that all positively selected B-cells stay in the GC for further rounds of mutation, division and selection. \\

Figure \ref{fig5_11} shows the dependence on $\overline a_s$ of the GC size and fitness, comparing $\mathcal M^+$ (left column) and $\mathcal M^-$ (right column). Indeed, for these models the GC  depends on the selection threshold, conversely to the previous case of positive and negative selection, and not only on the selection rate. The effects of $\overline a_s$ on the GC are perfectly simmetric: it is interesting to observe that when both selection mechanisms are coupled, then $\overline a_s$ does not affect the GC dynamics anymore, as shown for instance in Figure \ref{fig3_4} (a). Moreover, Figures \ref{fig5_11} (c,d) evidence the existence of a value of $\overline a_s$ that minimizes (resp. maximizes) the expected average affinity in the GC for $\mathcal M^+$ (resp. $\mathcal M^-$). In both cases this value is approximately $N/2$. This certainly depends on the transition probability matrix chosen for the mutational model, which converges to a binomial probability distribution over $\{0,\dots,N\}$. \\

\captionsetup[figure]{labelfont=bf}
\begin{figure}[ht!]
  \centering
  \begin{tabular}{cc}
  \subfloat[ ]{\label{fig11a}\includegraphics[scale=0.29]{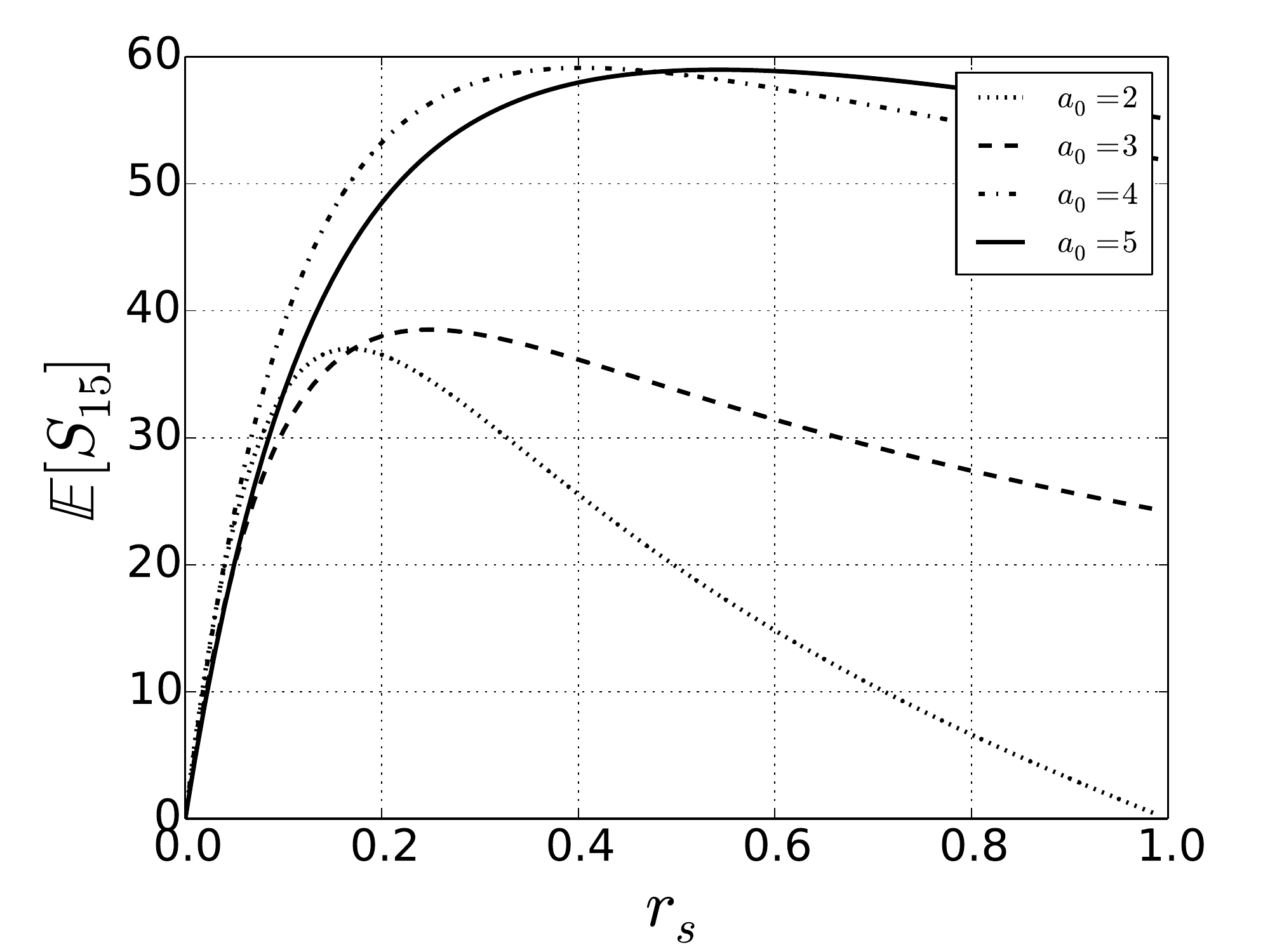}} &
  \subfloat[ ]{\label{fig11b}\includegraphics[scale=0.29]{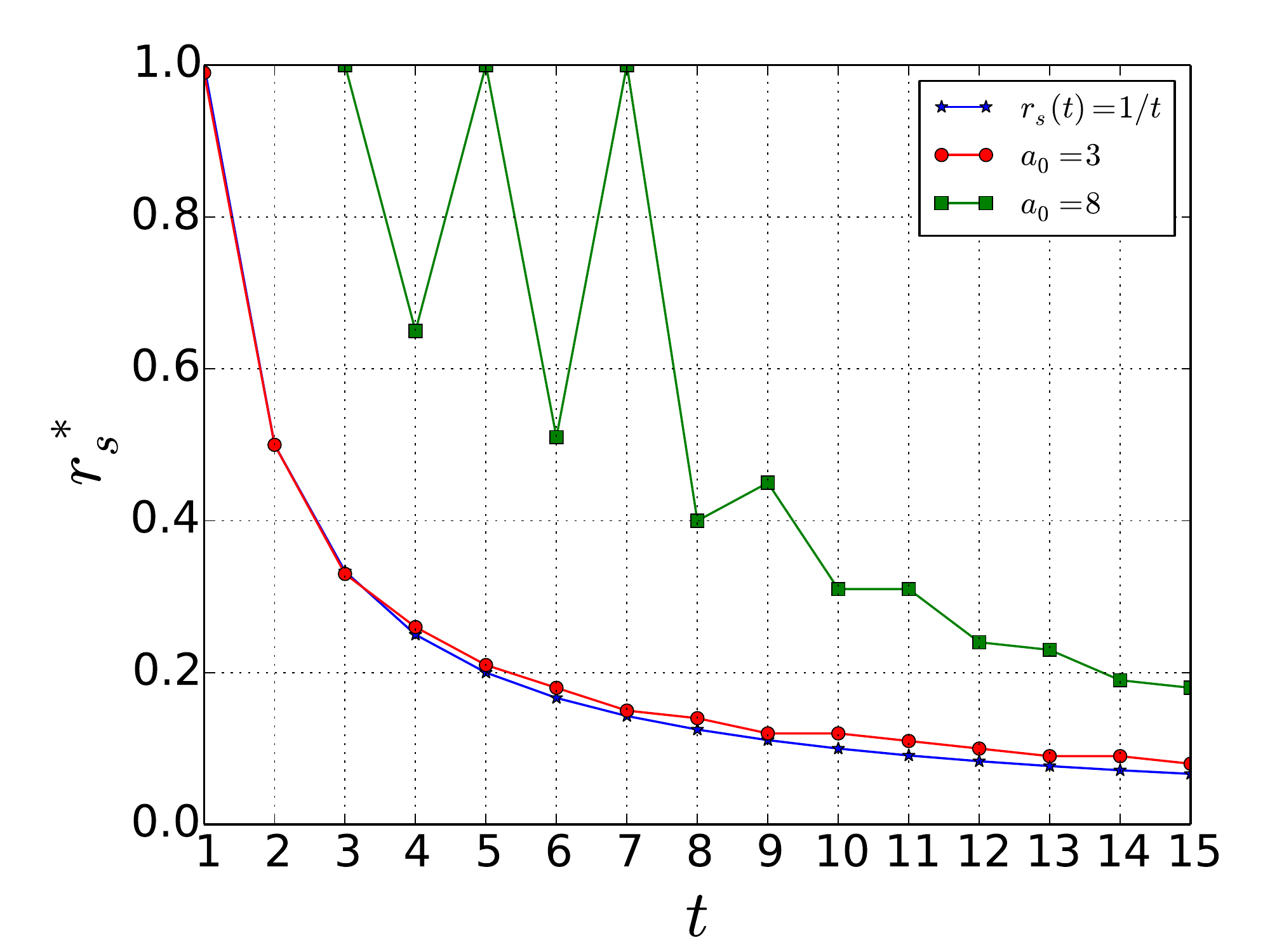}} 
  \end{tabular}
  \caption{Model of positive selection. (a) Expected number of selected B-cells for the time step $t=15$ for different values of $a_0$, depending on $r_s$. (b) Estimation of the optimal $r_s^\ast$ maximizing the expected number of selected B-cells for a given generation, comparing the model of positive selection for different values of $a_0$ and the model described in Section \ref{sec1} (we plot the exact value, $r_s(t)=1/t$, as obtained by Proposition \ref{cor:maxrs}). In (b), for simulations corresponding to the model of positive selection we set $\overline a_s=5$.}
  \label{fig11}
\end{figure}

The evolution of the selected pool for the model of positive selection have some important differences if compared to the model described in Section \ref{sec1}. For instance, it is not easy to identify an optimal value of $r_s$ which maximizes the expected number of selected B-cells at time $t$. Indeed it depends both on $a_0$ and $\overline a_s$: if $a_0\leq\overline a_s$ we find curves similar to those plotted in Figure \ref{fig9_10} (a), o\-ther\-wise Figure \ref{fig11} (a) shows a substantial different behavior. Indeed, if $a_0>\overline a_s$, choosing a big value for $r_s$ does not negatively affect the number of selected B-cells at time $t$. In this case, for the first time steps no (or a very few) B-cells will be positively selected, since they still need to improve their affinity to the target. Therefore, they stay in the GC and continue to proliferate for next generations. This fact is further underlined in Figure \ref{fig11} (b), where we estimate numerically the optimal $r_s^\ast$ which maximizes the expected number of selected B-cells at time $t$. Simulations show that for $a_0\leq\overline a_s$ the value of $r_s^\ast$ for the model of positive selection is really close to the one obtained by Proposition \ref{cor:maxrs}. On the other hand if we start from an initial affinity class $a_0>\overline a_s$ the result we obtain is substantially different from the previous one, especially for small $t$. Moreover we observe important oscillations, which are probably due to the mutational model, and to the fact that the total GC size is still small for small $t$, since the process starts from a single B-cell. Nevertheless, it seems that for $t$ big enough also in this case the value of $r_s^\ast$ tends to approach $1/t$.\\

\captionsetup[figure]{labelfont=bf}
\begin{figure}[ht!]
  \centering
  \begin{tabular}{cc}
  \subfloat[ ]{\label{fig16a}\includegraphics[scale=0.29]{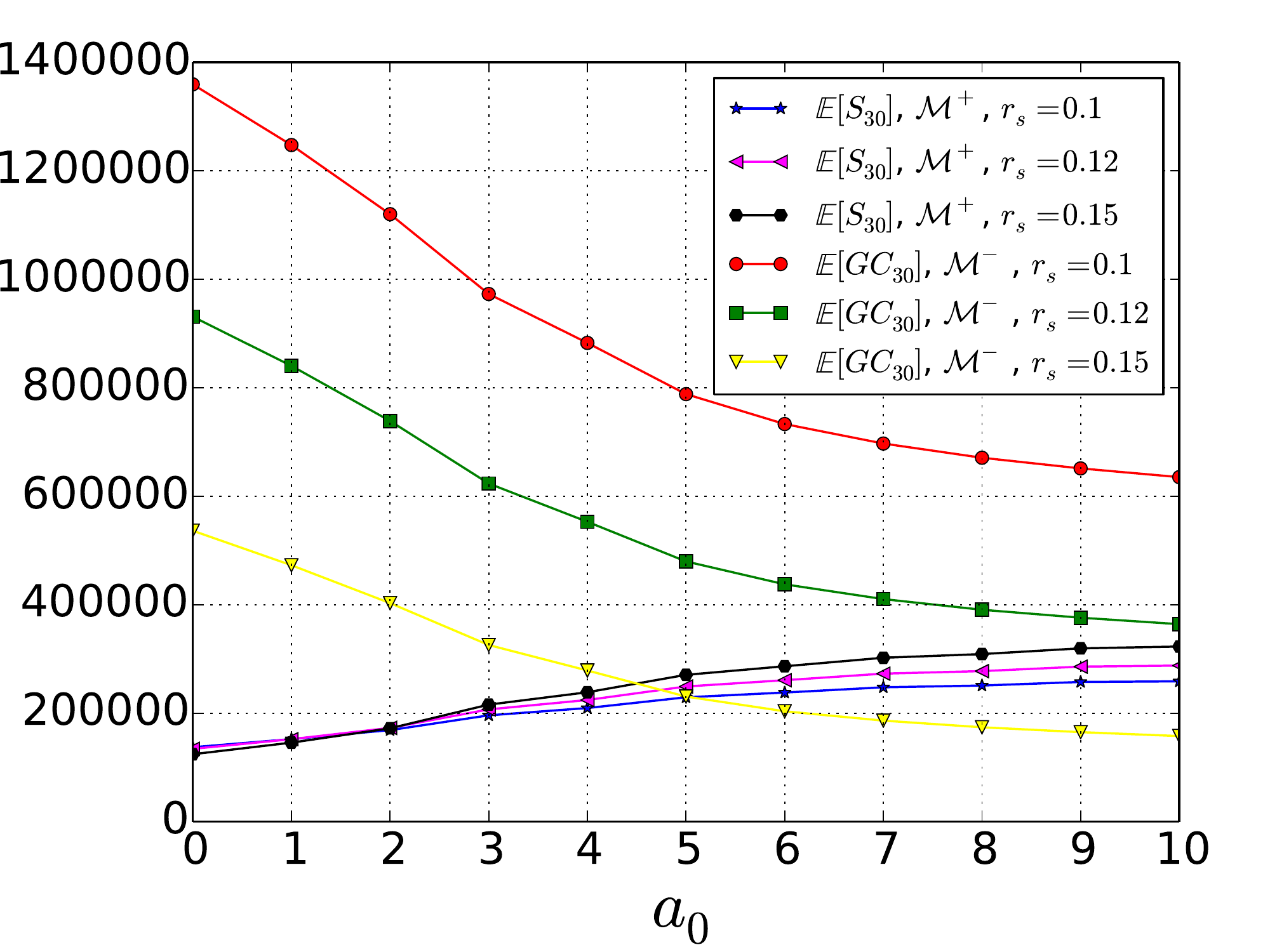}} &
  \subfloat[ ]{\label{fig16b}\includegraphics[scale=0.29]{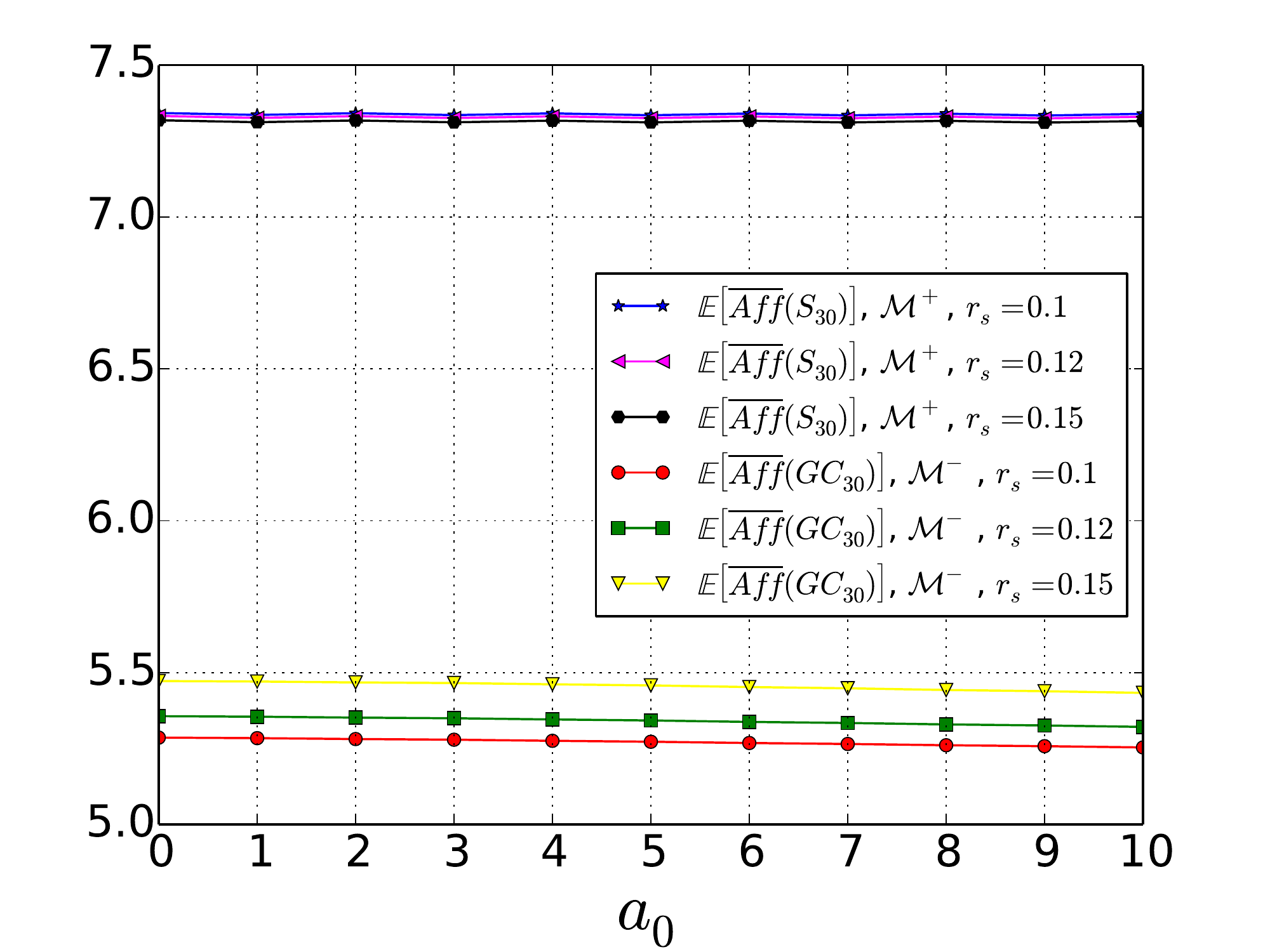}}
  \end{tabular}
  \caption{(a) Expected number of B-cells which have been selected until time $t=30$ for $\mathcal M^+$ compared to the expected size of the GC for $\mathcal M^-$ for different values of $r_s$. (b) Expected corresponding average affinity for the selected pool (case of positive selection) and the GC (case of negative selection). For some choice of the parameter $r_s$, the size of the selected pool for $\mathcal M^+$, and the GC for $\mathcal M^-$, are comparable. Nevertheless, the corresponding average affinities are significantly different.}
  \label{fig16}
\end{figure}

Since in the case of negative selection there is no selected pool, one can suppose that at a given time $t$ the process stops and all clones in the GC pool exit the GC as selected clones. Hence it can be interesting to compare the selected pool of the model of positive selection and the GC pool of the model of negative selection at time $t$. Clearly to make these two compartments comparable, the main parameters of both systems have to be opportunely chosen. In Figure \ref{fig16} we compare the size and average fitness of the selected pool for $\mathcal M^+$ and the GC for $\mathcal M^-$ at time $t=30$. We test different values of the parameter $r_s$. In particular, we observe that increasing $r_s$ the GC size for the model of negative selection decreases and its average fitness increases. For the parameter choices we made for these simulations, Figure \ref{fig16} (a) shows that the size of the GC for $\mathcal M^-$ is comparable to the size of the selected pool for $\mathcal M^+$ at time $t=30$ if, keeping all other parameters fixed, $r_s=0.15$ for $\mathcal M^-$. Nevertheless, this does not implies a comparable value for the average affinity: the clones of the selected pool for $\mathcal M^+$ have a significantly greater average affinity than those of the GC for $\mathcal M^-$. In order to increase the average fitness in the GC for the model of negative selection one has to consider greater values for the parameter $r_s$, but this affects the probability of extinction of the process.\\

We can expect this discrepancy between the average affinity for the selected pool for $\mathcal M^+$ and the one of the GC for $\mathcal M^-$. Indeed, in the first case we are looking to all those B-cells which have been positive selected, hence belong at most to the $\overline a_s^{\textrm{th}}$-affinity class. On the contrary in the case of $\mathcal M^-$, we consider the average affinity of all B-cells which are still alive in the GC at a given time step. Among these clones, if $r_s<1$, with positive probability there are also individuals with affinity smaller than the one required for escaping negative selection. They remain in the GC because they have not been submitted to selection. These B-cells make the average affinity decrease.
Of course $r_s$ is not the only parameter affecting the quantities plotted in Figure \ref{fig16}. In particular, one can observe that choosing a greater value for $\overline a_s$ also have a significant effect over the growth of both pools, as discussed in Remark \ref{remmaxlam}.

\section{Conclusions and perspectives}\label{sec5}

In this paper we formalize and analyze a mathematical model describing an evolutionary process with affinity-dependent selection. We use a multi-type GW process, obtaining a discrete-time probabilistic model, which includes division, mutation, death and selection. In the main model developed here, we choose a selection mechanism which acts both positively and negatively on individuals submitted to selection. This leads to build matrix $\mathcal M$, which contains the expectations of each type (Proposition \ref{propM}) and enables to describe the average behavior of all components of the process. Moreover, thanks to the spectral decomposition of $\mathcal M$ we were able to obtain explicitly some formulas giving the expected dynamics of all types. In addition, we exhibited an optimal value of the selection rate maximizing the expected number of selected clones for the $t^{\rm th}$-generation (Proposition \ref{cor:maxrs}). \\

This is one possible choice of the selection mechanism. From a ma\-the\-ma\-ti\-cal point of view, matrix $\mathcal M$ is particularly easy to manipulate, as we can obtain explicitly its spectrum. On the other hand, the positive and negative selection model leads, for example, to a selection threshold that does not have any impact on the evolution of the GC size. From a biological point of view this seems counterintuitive, since we could expect that the GC dynamics is sensible to the minimal fitness required for positive selection. Moreover, this process does not take into account any recycling mechanism, which has been confirmed by experiments \cite{victora2010germinal} and which improves  GCs' efficiency. In addition, we  considered  that only the selection mechanism is affinity dependent, while in the GC reaction  other mechanisms, such as the death and proliferation rate, may depend on fitness \cite{gitlin2014clonal,anderson2009taking}. Of course it is possible to define models with affinity-dependent division and death mechanisms with our formalism. This would clearly lead to a more complicated model, which can be at least studied numerically.\\

Mathematical tools used in Section \ref{sec1} can be applied to define and study other selection mechanisms. For instance in Section \ref{sec:ext} we propose two variants of the model analyzed in Section \ref{sec1}, in which selection acts only positively, resp. only negatively. This Section shows how our mathematical environment can be  modified to describe different selection mechanisms, which can be studied at least numerically. Moreover, it gives a deeper insight of the previous model of positive and negative selection, by highlighting the effects of each selection mechanism individually, when they are not coupled. \\

From a biological viewpoint there exist many possibilities to improve the models proposed in this paper. First of all it is extremely important to fix the system parameters, which have to be consistent with the real biological process. The choice of $N$ defines the number of affinity level with respect to  a given antigen. This value can be interpreted in different ways. On the one hand it can correspond to the number of key mutations observed during the process of Antigen Affinity Maturation, hence be even smaller than 10.  On the other hand, each mutational event implies a change in the B-cell affinity, slight or not if it is a key mutation. In this case the affinity can be modeled as a continuous function, hence $N$ corresponds to a possible discretization \cite{weiser2011affinity,xu2015key}. To this choice corresponds an appropriate choice of the transition probability matrix defining the mutational model over the affinity classes, $\mathcal Q_N$. In most numerical simulations we set $N=10$, which is a sensible value since experimentalists observe that high-affinity B-cells differ in their BCR coding gene by about 9 mutations from germline genes \cite{DIbePKMai,zhang2010optimality}.
Nevertheless all ma\-the\-ma\-ti\-cal results are independent from this choice and hold true for all $N\geq1$. The selection, division and death rates have also an important impact in the GC and selected pool dynamics: in the simulations we set them in order to be in a case of explosion of the GC hence appreciate the effects of all parameters over the main quantities, but they are not biologically justified. 
For instance, the typical proliferation rate of a B-cell has been estimated between 2 and 4 per day and in the literature we found B-cell death rates of the order of 0.5-0.8 per day \cite{meyer2006analysis,zhang2010optimality,kecsmir1999mathematical}. Hence, if we suppose that a single time step corresponds \emph{e.g.} to 6 hours, a consistent proliferation rate would be $r_{div}\simeq0.75$, while the death rate $r_d$ should be around $0.175$. Since over a 6 hours period about 50$\%$ of B-cells transit from the DZ to the LZ, where they compete for positive selection signaling \cite{bannard2013germinal,victora2014snapshot}, we should choose $r_s\leq0.5$. It could be further characterized taking into account its tightly relation with the time of GC peak, as highlighted in Section \ref{sec:rsbest}.\\

In Section \ref{sec:rsbest} we have explicitly determined the optimal value of the selection rate maximizing the production of output cells at time $t$ for the main model of positive and negative selection. It is equal to $1/t$ independently from all other parameters. Moreover, numerical estimations for the model of positive selection (Section \ref{sec5:2}) suggest that also in this case there exists an optimal value of $r_s(t)$, which tends to $1/t$ at least for $t$ big enough. One has to interpret this result as the ideal optimal strength of the selection pressure to obtain a peak of the GC production of output cells at a given time step. For example, let us suppose again that a time step corresponds to 6 hours. The peak of the GC reaction has been measured to be close to day 12 \cite{wollenberg2011regulation}, \emph{i.e.} after $\sim48$ maturation cycles in our model: for the kind of models we built and analyzed in this paper, a constant selection pressure $r_s$ of $1/48\simeq0.02$ assures that the production of plasma and memory B-cells at the GC peak is maximized. Note that with the parameter choice $r_d=0.175$, $r_{div}=0.75$ and $r_s=0.02$, the extinction probability of the GC is $\simeq0.3^{z_0}$, $z_0$ being the number of initial seed cells. 
Since the extinction probability is strictly smaller than 1, such a GC will explode with high probability and will be able to assure an intense and efficient immune response. \\

\captionsetup[figure]{labelfont=bf}
\begin{figure}[ht!]
  \centering
  \subfloat[ ]{\label{fig17a}\includegraphics[scale=0.29]{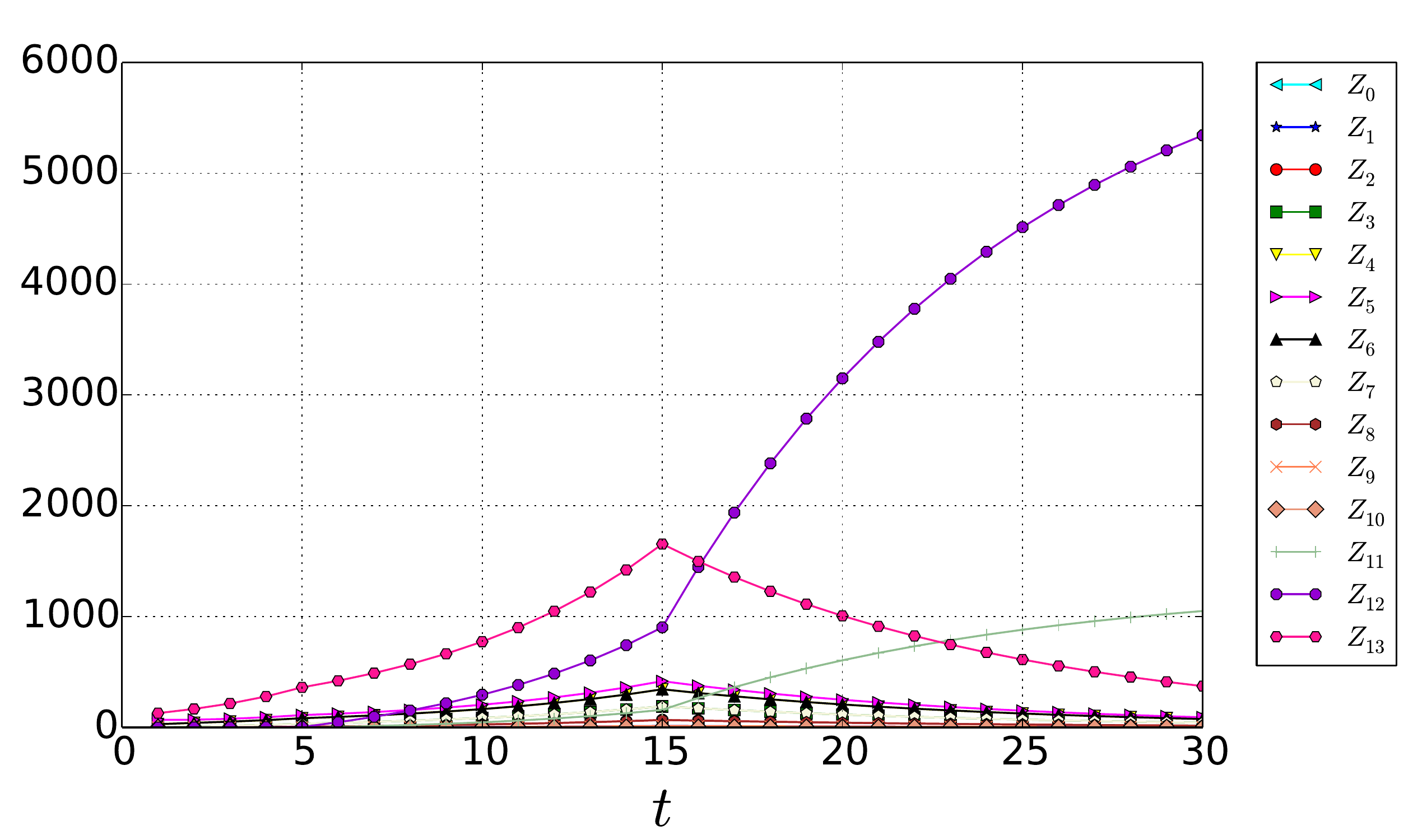}}~
  \subfloat[ ]{\label{fig17b}\includegraphics[scale=0.29]{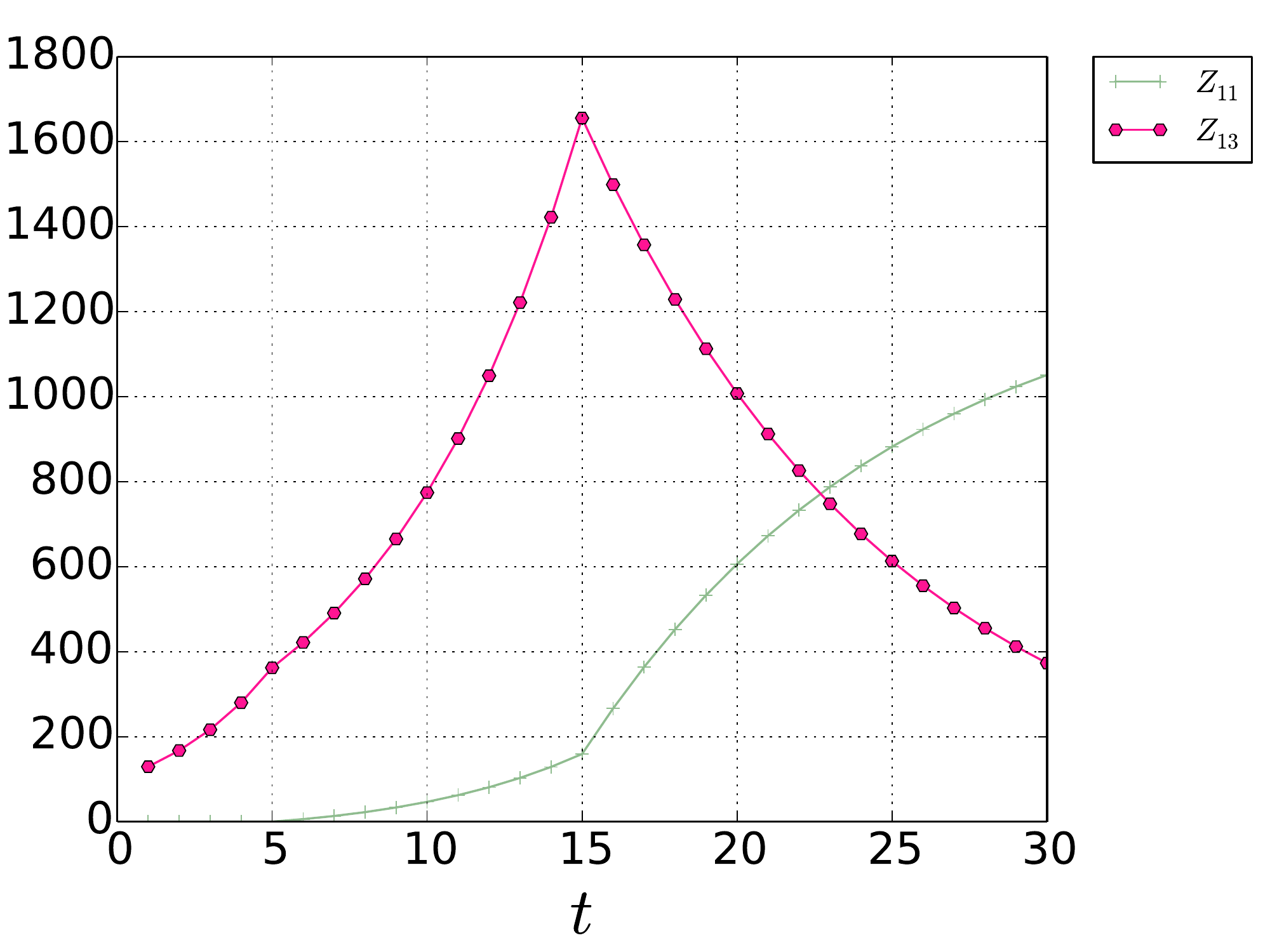}}
  \caption{(a) Evolution during time of the expected value of all types for the model of positive and negative selection, with $r_s$ varying during time and $N=10$. In particular we set $r_s=0$ until $t= 5$, $r_s=0.1$ from $t=6$ to $t=15$ and $r_s=0.3$ from $t=16$ to $t=30$. $Z_{13}$ denotes the total size of the GC (\emph{i.e.} $\sum_{k=0}^NZ_k$), and we recall that $Z_{11}$ corresponds to selected B-cells and $Z_{12}$ to dead B-cells. We set $r_{div}=0.3$, $r_d=0.005$ and $z_0=100$ initial naive B-cells. All initial B-cells belong to $a_0=5$, and the selection threshold is $\overline a_s=3$. (b) Evolution during time of the expected total size of the GC and the selected pool respectively, for the same set of parameters as in Figure \ref{fig17} (a). }
  \label{fig17}
\end{figure}

In our models  the selection pressure is constant. 
Since the optimal selection rate above depends on time, this suggests to go further in this direction.
Moreover, a time-dependent selection pressure  would allow to take into account, for instance, the early GC phase in which  simple clonal expansion of B-cells with no selection occurs \cite{de2015dynamics}. 
The hypothesis of a selection pressure changing over time can be easily integrated in our model. 
Indeed let us suppose that a selection rate $r_{s,1}$ until time $t_1$ and $r_{s,2}$ for all $t>t_1$ are fixed. 
Starting from the initial condition $\vec i$ the expectations of each type at time $t$ are given by $(\vec i\mathcal M_{r_{s,1}}^t)$ if $t\leq t_1$ and $(\vec i\mathcal M_{r_{s,1}}^{t_1}\mathcal M_{r_{s,2}}^{t-t_1})$ if $t>t_1$, where $\mathcal M_{r_{s,i}}$ is the matrix containing the expectations of each type for an evolutionary process with constant selection rate $r_{s,i}$, $i=1,2$. In Figure \ref{fig17} we plot the expected evolution during time of all types considering an increasing selection rate. We evaluate the expectations of all types following a process with positive and negative selection. We set $r_s=0$ until $t= 5$, $r_s=0.1$ from $t=6$ to $t=15$ and $r_s=0.3$ for $t>15$. 
Numerical simulations  show that a time dependent  selection rate  allows initial explosion of the GC, and then progressive extinction, 
while when  parameters are fixed,  a GW process  gives only  rise either to explosion or to extinction, as shown above. 
The regulation and termination of  the GC reaction has not yet been fully understood. 
In the literature, an increasing differentiation rate of GC B-cells is thought to be a good explanation \cite{moreira2006modelling}, here we show
that other reasons could be of importance as well.
Similarly, we can let  other parameters vary for fixed time intervals, as well as decide to alternatively switch on and off the mutational mechanism, as already proposed in \cite{ASPerGWeis}. This can be obtained by alternatively use the identity matrix in place of $\mathcal Q_N$. \\

Applications of the models presented here to real biological problems and data should be further investigated. We propose here some contexts for which we believe that our kind of modeling approach could be employed to address biologically relevant questions. \\
Even if it is still extremely hard to have precise experimental information about the evolution of Antibody Affinity Ma\-tu\-ra\-tion inside GCs, new refined techniques start to be available to measure clonal diversity in GCs. As an example, in \cite{tas2016visualizing} the authors combine multiphoton microscopy and sequencing to understand how different clonal diversification patterns can lead to efficient affinity maturation. The models we propose could be used to infer which are reasonable mutational transitional probability matrices and selection mechanisms/pressure to obtain such different scenario and infer if the tendency of GC to go or not through homogenizing selection is solely due to the hazard or if this is dependent on the kind of antigenic challenge and/or some specific characteristics of the host. If this is the case, these results could be particularly relevant \emph{e.g.} in the context of vaccination design, where we are interested in find new way to improve the quality of the immune response after vaccination challenge. \\
Another potential interesting application field is the study of some diseases entailing a dysfunction of the immune system, such as in particular Chronic Lymphocytic Leukemia (CLL), derived from antigen-experienced B-cells that differ in the level of mutations in their receptors \cite{chiorazzi2005chronic}. This is the commonest form of leukemia in the Western world \cite{eichhorst2015chronic}. In CLL, leukemia B-cells can mature partially but not completely, are unable to opportunely undergo mutations in GCs, and survive longer than normal cells, crowding out healthy B-cells. Prognosis varies depending on the ability of host B-cells to mutate their antibody gene variable region. Even if major progresses have been made in the identification of molecular and cellular markers   predicting the expansion of this disease  in patients, the pathology remains incurable \cite{dighiero2008chronic,eichhorst2015chronic}. Our mo\-deling approach could be employed to understand how an ``healthy'' mutational matrix is modified in patients affected by CLL, and if other mechanisms could contribute to get the prognosis worse. This could eventually provide suggestions about the causes that lead to CLL, and motivation for further research on possible treatments.

%
%


\bibliographystyle{spmpsci}      

\section*{Appendix}

\appendix

\section{Few reminders of classical results on GW processes}\label{app:classicGW}

We recall here some classical results about GW processes we employed to derive Proposition \ref{prop:eta} (Section \ref{sec2}). For further details the reader can refer to \cite{harris2002theory}.

\begin{definition}
Let $X$ be an integer valued rv, $p_k:=\Pro(X=k)$ for all $k\geq0$. Its probability generating function (pgf) is given by:
\begin{equation*}
F_X(s)=\sum_{k=0}^{+\infty} p_k s^k
\end{equation*}
\end{definition}

$F_X$ is a convex monotonically increasing function over $[0,1]$, and $F_X(1)=1$. If $p_0\neq0$ and $p_0+p_1<1$ then $F$ is a strictly increasing function.

\begin{definition}
Given $F$, the pgf of a rv $X$, the iterates of $F$ are given by:
\begin{equation*}
\left.\begin{array}{l}
F_0(s)=s \\
F_1(s)=F(s) \\
F_t(s)=F(F_{t-1}(s))\textrm{ for $t\geq2$}
\end{array}\right.
\end{equation*} 
\end{definition}

\begin{proposition}\label{pro:FEV}
\begin{description}
\item[]
\item[(i)] If $\E(X)$ exists (respectively $\mathbb V(X)$), then $\E(X)=F'_X(1)$ (respectively $\mathbb V(X)=F''_X(1)-\left(\E(X)\right)^2+\E(X)$).
\item[(ii)]If $X$ and $Y$ are two integer valued independent rvs, then $X+Y$ is still an integer valued rv and its pgf is given by $F_{X+Y}=F_XF_Y$.
\end{description}
\end{proposition}

\begin{definition}
We denote by $\eta$ the extinction probability of the process $(Z_t)_{t\in\mathbb N}$:
\begin{equation*}
\eta:=\lim_{t\to\infty}F_t(0)
\end{equation*}
\end{definition}

\begin{theorem}\label{thm:etaE}
\begin{description}
\item[]
\item[(i)] The pgf of $Z_t^{(z_0)}$, $t\;\in\;\mathbb N$, which represents the population size of the $t^{\textrm{th}}$-generation starting from $z_0\geq1$ seed cells, is $F_t^{(z_0)}=(F_t)^{z_0}$, $F_t$  being the $t^{\textrm{th}}$-iterate of $F$ (Equation \eqref{eq:gen}).
\item[(ii)] The expected size of the GC at time $t$ and starting from $z_0$ B-cells is given by:
\begin{equation}\label{eq:expZtz0}
\E(Z_t^{(z_0)})=z_0\left(\E(Z_t)\right)=z_0\left(\E(Z_1)\right)^t~,
\end{equation}
\item[(iii)]$\eta$ is the smallest fixed point of the generating function $F$, i.e. $\eta$ is the smallest $s$ s.t. $F(s)=s$.
\item[(iv)] If $\E(Z_1)=:m$ is finite, then:
\begin{itemize}
\item if $m\leq 1$ then $F$ has only 1 as fixed point and consequently $\eta=1$;
\item if $m>1$ then $F$ as exactly a fixed point on $[0,1[$ and then $\eta<1$.
\end{itemize}
\item[(v)] Denoted by $\eta_{z_0}$ the probability of extinction of $(Z_t^{(z_0)})$, one has:
\begin{equation*}
\eta_{z_0}=\eta^{z_0}
\end{equation*}
where $\eta$ is given by (iii).
\end{description}
\end{theorem}

Proposition \ref{prop:eta} of Section \ref{sec2} follows by applying Theorem \ref{thm:etaE} and Equation \eqref{eq:proZt}.

\section{Proof of Proposition \ref{propM}}\label{app:proofPropM}

For all $j\,\in\,\{0,\dots,N+2\}$ the generating function of $Z_j$ gives the number of offsprings of each type that a type $j$ particle can produce. It is defined as follows:
\begin{equation}
f^{(j)}(s_0,\dots,s_{N+2})=\displaystyle\sum_{k_0,\dots,k_{N+2}\geq 0}p^{(j)}(k_0,\dots,k_{N+2})s_0^{k_0}\cdots s_{N+2}^{k_{N+2}},
\end{equation}
\begin{equation*}
0\leq s_{\alpha}\leq1\textrm{ for all }\alpha\,\in\,\{0,\dots,N+2\}
\end{equation*}
where $p^{(j)}(k_0,\dots,k_{N+2})$ is the probability that a type $j$ cell produces $k_0$ cells of type $0$, $k_1$ of type $1$, $\dots$, $k_{N+2}$ of type $N+2$ for the next generation.\\
We denote:
\begin{itemize}
\item $\mathbf p(\mathbf k)=(p^{(0)}(\mathbf k),\dots,p^{(N+2)}(\mathbf k))$, for $\mathbf k=(k_0,\dots,k_{N+2}) \,\in\,\mathbb Z_+^{N+3}$
\item $\mathbf f(\mathbf s)=(f^{(1)}(\mathbf s),\dots,f^{(N+1)}(\mathbf s))$, for $\mathbf s=(s_0,\dots,s_{N+2})\,\in\,\mathcal C^{N+3}:=[0,1]^{N+3}$
\end{itemize}
Then the probability generating function of $\mathbf Z_1$ is given by:
\begin{equation}
\mathbf f(\mathbf s)=\displaystyle\sum_{\mathbf k\in\mathbb Z_+^{N+3}}\mathbf p(\mathbf k)\mathbf s^{\mathbf k}\textrm{, $\mathbf s\,\in\,\mathcal C^{N+3}$}
\end{equation}
Again, the generating function of $\mathbf Z_t$, $\mathbf f_t(\mathbf s)$, is obtained as the $t^{\textrm{th}}$-iterate of $\mathbf f$, and it holds true that:
\begin{equation*}
\mathbf f_{t+r}(\mathbf s)=\mathbf f_t[\mathbf f_r(\mathbf s)]\textrm{, $\mathbf s\,\in\,\mathcal C^{N+3}$.}
\end{equation*}

Let $m_{ij}:=\E[Z_{1,j}^{(i)}]$ the expected number of offspring of type $j$ of a cell of type $i$ in one generation. We collect all $m_{ij}$ in a matrix, $\mathcal M=(m_{ij})_{0\leq i,j\leq N+2}$. We have \cite{athreya2012branching}:
\begin{equation*}
m_{ij}=\frac{\partial f^{(i)}}{\partial s_j}(\boldsymbol 1)
\end{equation*}
and:
\begin{equation}
\E[Z_{t,j}^{(i)}]=\frac{\partial f_t^{(i)}}{\partial s_j}(\boldsymbol 1)
\end{equation}
Finally:
\begin{equation}\label{eq:EZn}
\E[\mathbf Z_t^{(\mathbf i)}]=\mathbf i\mathcal M^t
\end{equation}

One can explicitly derive the elements of matrix $\mathcal M$ for the process described in Definition \ref{def:Ztmulti}.

\begin{proposition*}
$\mathcal M$ is a $(N+3)\times(N+3)$ matrix defined as a block matrix:
\begin{equation*}
\mathcal M=\left(\begin{array}{cc}
\mathcal M_1 & \mathcal M_2 \\
\boldsymbol 0_{2\times (N+1)} & \mathcal I_2
\end{array}\right)
\end{equation*}
Where:
\begin{itemize}
\item $\boldsymbol 0_{2\times (N+1)}$ is a $2\times (N+1)$ matrix with all entries 0;
\item $\mathcal I_n$ is the identity matrix of size $n$;
\item $\mathcal M_1=2(1-r_d)r_{div}(1-r_s)\mathcal Q_N+(1-r_d)(1-r_{div})(1-r_s)\mathcal I_{N+1}$
\item $\mathcal M_2=(m_{2,ij})$ is a $(N+1)\times 2$ matrix where for all $i\,\in\,\{0,\dots,N\}$:
\begin{itemize}
\item if $i\leq\overline a_s$: \\
$m_{2,i1}=(1-r_d)(1-r_{div})r_s+2(1-r_d)r_{div}r_s\displaystyle\sum_{j=0}^{\overline a_s}q_{ij}$, \\
$m_{2,i2}=r_d+2(1-r_d)r_{div}r_s\displaystyle\sum_{j=\overline a_s+1}^{N}q_{ij}$
\item if $i>\overline a_s$: \\
$m_{2,i1}=2(1-r_d)r_{div}r_s\displaystyle\sum_{j=0}^{\overline a_s}q_{ij}$, \\
$m_{2,i2}=r_d+(1-r_d)(1-r_{div})r_s+2(1-r_d)r_{div}r_s\displaystyle\sum_{j=\overline a_s+1}^{N}q_{ij}$
\end{itemize}
\end{itemize}
\end{proposition*}

\proof
One has to compute all $f^{(i)}(\mathbf s)$ for $i=0,\dots, N+2$, 
which depend on $r_d$, $r_{div}$, $r_s$, $\overline a_s$ and the elements of $\mathcal Q_N$. 
First, the elements of the $(N+2)^{\textrm{th}}$ and $(N+3)^{\textrm{th}}$-lines are obviously determined: all selected (resp. dead) cells remain selected (resp. dead) for next generations, as they can not give rise to any other cell type offspring (we do not take into account here any type of recycling mechanism). Let $i\;\in\;\{0,\dots,N\}$ be a fixed index: we evaluate $m_{ij}$ for all $j\;\in\;\{0,\dots,N+2\}$. 
The first step is to determine the value of $p^{(i)}(\mathbf k)$ for $\mathbf k=(k_0,\dots,k_{N+2}) \,\in\,\mathbb Z_+^{N+3}$. There exists only a few cases in which $p^{(i)}(\mathbf k)\neq0$, which can be explicitly evaluated:
\begin{itemize}
\item $p^{(i)}(0,\dots,0,1)=\begin{cases}
 r_d & \text{if}\;  i \leq \overline a_s \\
 r_d+(1-r_d)(1-r_{div})r_s & \text{otherwise}
\end{cases}$
\item $p^{(i)}(0,\dots,0,1,0)=\begin{cases}
 (1-r_d)(1-r_{div})r_s & \text{if}\;  i \leq \overline a_s \\
 0 & \text{otherwise}
\end{cases}$
\item $p^{(i)}(0,\dots,0,\underset{{\color{gray50} i}}{1},0,\dots,0,0)=(1-r_d)(1-r_{div})(1-r_s)$
\item $p^{(i)}(0,\dots,0,2)=(1-r_d)r_{div}r_s^2\displaystyle\sum_{j_1=\overline a_s+1}^N q_{ij_1}\sum_{j_2=\overline a_s+1}^N q_{ij_2}$
\item $p^{(i)}(0,\dots,0,2,0)=(1-r_d)r_{div}r_s^2\displaystyle\sum_{j_1=0}^{\overline a_s} q_{ij_1}\sum_{j_2=0}^{\overline a_s} q_{ij_2}$
\item $p^{(i)}(0,\dots,0,1,1)=2(1-r_d)r_{div}r_s^2\displaystyle\sum_{j_1=0}^{\overline a_s} q_{ij_1}\sum_{j_2=\overline a_s+1}^N q_{ij_2}$
\item For all $j_1<j_2\;\in\;\{0,\dots,N\}$:
\begin{itemize}
\item $p^{(i)}(0,\dots,0,\underset{{\color{gray50} j_1}}{2},0,\dots,0,0)=(1-r_d)r_{div}(1-r_s)^2q_{ij_1}^2$
\item $p^{(i)}(0,\dots,0,\underset{{\color{gray50} j_1}}{1},0,\dots,0,\underset{{\color{gray50} j_2}}{1},0,\dots,0,0)=2(1-r_d)r_{div}(1-r_s)^2q_{ij_1}q_{ij_2}$
\item $p^{(i)}(0,\dots,0,\underset{{\color{gray50} j_1}}{1},0,\dots,0,1)=2(1-r_d)r_{div}r_s(1-r_s)q_{ij_1}\displaystyle\sum_{j_2=\overline a_s+1}^N q_{ij_2}$
\item $p^{(i)}(0,\dots,0,\underset{{\color{gray50} j_1}}{1},0,\dots,0,1,0)=2(1-r_d)r_{div}r_s(1-r_s)q_{ij_1}\displaystyle\sum_{j_2=0}^{\overline a_s} q_{ij_2}$
\end{itemize}
\item $p^{(i)}(\mathbf k)=0$ otherwise
\end{itemize}
We can therefore evaluate $f^{(i)}(\mathbf s)$, with $\mathbf s=(s_0,\dots,s_{N+2})\,\in\,\mathcal C^{N+3}$. \\

For all $i\leq\overline a_s$:
\begin{equation} \label{eq:fis}
 \hspace{-0.15cm} \begin{aligned}
& f^{(i)} (\mathbf s)  =  r_ds_{N+2}+(1-r_d)(1-r_{div})r_ss_{N+1} + (1-r_d)(1-r_{div})(1-r_s)s_i \\
& +  (1-r_d)r_{div}r_s^2\left(\displaystyle\sum_{j_1=\overline a_s+1}^N q_{ij_1}\sum_{j_2=\overline a_s+1}^N q_{ij_2}s_{N+2}^2\right.\\
& \left. +\sum_{j_1=0}^{\overline a_s} q_{ij_1}\sum_{j_2=0}^{\overline a_s} q_{ij_2}s_{N+1}^2+2\sum_{j_1=0}^{\overline a_s} q_{ij_1}\sum_{j_2=\overline a_s+1}^N q_{ij_2}s_{N+1}s_{N+2}\right)   \\  \\
& +  (1-r_d)r_{div}(1-r_s)^2\left(\displaystyle\sum_{j_1=0}^{N}q_{ij_1}^2s_{j_1}^2+2\sum_{j_1=0}^{N}q_{ij_1}\sum_{j_2< j_1=0}^N q_{ij_2}s_{j_1}s_{j_2}\right)  \\  \\
& +  2(1-r_d)r_{div}r_s(1-r_s)\displaystyle\sum_{j_1=0}^{N}q_{ij_1}\left(\displaystyle\sum_{j_2=\overline a_s+1}^N q_{ij_2}s_{N+2}+\sum_{j_2=0}^{\overline a_s} q_{ij_2}s_{N+1}\right)s_{j_1}
\end{aligned}
\end{equation}

If $i>\overline a_s$ then $f^{(i)}(\mathbf s)$ is the same except for the first line, which becomes:
\begin{equation*}
(r_d+(1-r_d)(1-r_{div})r_s)s_{N+2} + (1-r_d)(1-r_{div})(1-r_s)s_i
\end{equation*}
The values of each $m_{ij}$ are now obtained by evaluating all partial derivatives of $f^{(i)}(\mathbf s)$ in $\boldsymbol 1$, keeping in mind that for all $i\;\in\;\{0,\dots,N\}$, $\sum_{j=0}^N q_{ij}=1$.
\qed
\endproof

\section{Deriving the extinction probability of the GC from the multi-type GW process (Section \ref{sec3})}\label{app:extinctionProbmulti}

Let us recall some results about the extinction probability for multi-type GW processes \cite{athreya2012branching}.

\begin{definition}\label{def:extmulti}
Let $q^{(i)}$ be the probability of eventual extinction of the process, when it starts from a single type $i$ cell. 
As above bold symbols denote vectors {\em i.e.}
$\mathbf q :=(q^{(0)},\dots,q^{(N+2)})\geq 0$.
\end{definition}

\begin{definition}
We say that $(\mathbf Z_t)$ is singular if each particle has exactly one offspring, which implies that the branching process becomes a simple MC. 
\end{definition}

\begin{definition}\label{def:pos}
Matrix $\mathcal M$ is said to be strictly positive if it has non-negative entries and there exists a $t$ s.t. $\left(\mathcal M^t\right)_{ij}>0$ for all $i$, $j$. $(\mathbf Z_t)$ is called positive regular iff $\mathcal M$ is strictly positive.
\end{definition}

\begin{notation}
Let $\mathbf u$, $\mathbf v\;\in\;\mathbb R^n$. We say that $\mathbf u\leq \mathbf v$ if $u_i\leq v_i$ for all $i\;\in\:\{1,\dots,n\}$. Moreover, we say that $\mathbf u< \mathbf v$ if $\mathbf u\leq \mathbf v$ and $\mathbf u\neq \mathbf v$.
\end{notation}

\begin{theorem}\label{thm:extmultit}
Let $(\mathbf Z_t)$ be non singular and strictly positive. Let $\rho$ be the maximal eigenvalue of $\mathcal M$. The following three results hold:
\begin{enumerate}
\item If $\rho<1$ (subcritical case) or $\rho=1$ (critical case) then $\mathbf q=\boldsymbol 1$. Otherwise, if $\rho>1$ (supercritical case), then $\mathbf q<\boldsymbol 1$.
\item $\displaystyle\lim_{t\to\infty}\mathbf f_t(\mathbf s)=\mathbf q$, for all $\mathbf s\,\in\,\mathcal C^{N+3}$.
\item $\mathbf q$ is the only solution of $\mathbf f(\mathbf s)=\mathbf s$ in $\mathcal C^{N+3}$.
\end{enumerate}
\end{theorem}

The spectrum of matrix $\mathcal M$ defined in Definition \ref{propM}  (and recalled in Appendix \ref{app:proofPropM}) is obtained as follows:

\begin{proposition}
Let $\mathcal M$ be defined as a block matrix as in  Proposition \ref{propM}. Let $\lambda_{\mathcal M,i}$ be its $i^{\textrm{th}}$-eigenvalue. The spectrum of $\mathcal M$ is given by:
\begin{itemize}
\item For all $i\;\in\;\{0,\dots,N\}$, $\lambda_{\mathcal M,i}=(1-r_d)(1-r_s)(1+r_{div}(2\lambda_{i}-1))$, where $\lambda_{i}$ is the $i^{\textrm{th}}$-eigenvalue of matrix $\mathcal Q_N$.
\item whereas $\lambda_{\mathcal M,N+1}=1$ with multiplicity 2.
\end{itemize}
\end{proposition}

\proof
As $\mathcal M$ is a block matrix with the lower left block composed of zeros, then $Spec(\mathcal M)=Spec(\mathcal M_1)\cup Spec(\mathcal I_{2})$. The result follows.
\qed
\endproof

Therefore we obtain the same condition as in Proposition \ref{prop:eta} for the extinction probability in the GC:

\begin{proposition}\label{prop:extGCpm}
Let $\mathbf q$ be the extinction probability for the process $(\mathbf Z_t)$ defined in Definition \ref{def:Ztmulti} and restricted to the first $N+1$ components (\emph{i.e.} we refer only to matrix $\mathcal M_1$, which defines the expectations of GC B-cells). Therefore:
\begin{itemize}
\item if $r_s\geq1-\displaystyle\frac{1}{(1-r_d)(1+r_{div})}$, then $\mathbf q=\boldsymbol 1$
\item otherwise $\mathbf q<\boldsymbol 1$ is the smallest fixed point of $\mathbf f(\mathbf s)$ in $\mathcal C^{N+3}$.
\end{itemize}
\end{proposition}

\proof
$\mathcal Q_N$ is a stochastic matrix, therefore its largest eigenvalue is 1. The corresponding eigenvalue of matrix $\mathcal M_1$ is: $\lambda_{\mathcal M_1,1}=(1-r_d)(1-r_s)(1+r_{div})$. The proposition is proved by observing that $\lambda_{\mathcal M_1,1}\leq1\Leftrightarrow r_s\geq1-\displaystyle\frac{1}{(1-r_d)(1+r_{div})}$ and applying Theorem \ref{thm:extmultit} (note that $\mathcal M_1$ is positive regular: this is not the case for matrix $\mathcal M$).
\qed
\endproof

\section{Expected size of the GC derived from the multi-type GW process (Section \ref{sec3})}\label{app:proofeq:sizeGC}

\begin{proposition*}
Let $\mathbf i$ be the initial state, $z_0:=|\mathbf i|$ its 1-norm ($|\mathbf i|:=\sum_{j=0}^{N+2}\mathbf i_j$). The expected size of the GC at time $t$:
\begin{equation*}
\displaystyle\sum_{k=0}^N(\mathbf i\mathcal M^t)_k=|\mathbf i|\left((1-r_d)(1+r_{div})(1-r_s)\right)^{t}
\end{equation*}
\end{proposition*}

\proof
For the sake of simplicity, let us suppose that the process starts from a single B-cell belonging to the affinity class $a_0=i$ with respect to  the target trait. We do not need to specify the transition probability matrix used to define the mutational model allowed. \\

We recall the expression of $\mathcal M^t$ obtained by iteration:
\begin{equation*}
{\mathcal M}^t=\left(\begin{array}{cc}
\mathcal M_1^t & \displaystyle\sum_{k=0}^{t-1}\mathcal M_1^k\mathcal M_2 \\ \\
\boldsymbol 0_{2\times (N+1)} & \mathcal I_2
\end{array}\right)
\end{equation*}

Therefore we can claim that $(\mathbf i\mathcal M^t)_k$ corresponds to the $k^{\textrm{th}}$-component of the $i^{\textrm{th}}$-row of matrix $\mathcal M_1^t=(2(1-r_d)r_{div}(1-r_s)\mathcal Q_N+(1-r_d)(1-r_{div})(1-r_s)\mathcal I_{N+1})^t$, where $\mathcal Q_N$ is a stochastic matrix. Matrices $\mathcal A:=2(1-r_d)r_{div}(1-r_s)\mathcal Q_N$ and $\mathcal B:=(1-r_d)(1-r_{div})(1-r_s)\mathcal I_{N+1}$ clearly commute, therefore:
\begin{equation}
\left(\mathcal A+\mathcal B\right)^t=\sum_{j=0}^tC_t^j\mathcal A^{t-j}\mathcal B^j
\end{equation}

For all $j$, $0\leq j\leq t$:
\begin{eqnarray*}
\mathcal A^{t-j}\mathcal B^j & = & 2^{t-j}(1-r_d)^{t-j}r_{div}^{t-j}(1-r_s)^{t-j}(1-r_d)^j(1-r_{div})^j(1-r_s)^j\mathcal Q_N^{t-j} \\
& = & (1-r_d)^t(1-r_s)^t(2r_{div})^{t-j}(1-r_{div})^j\mathcal Q_N^{t-j}
\end{eqnarray*}
Hence:
\begin{equation*}
\left(\mathcal A+\mathcal B\right)^t=(1-r_d)^t(1-r_s)^t\sum_{j=0}^tC_t^j(2r_{div})^{t-j}(1-r_{div})^j\mathcal Q_N^{t-j}
\end{equation*}
And consequently:
\begin{eqnarray*}
\sum_{k=0}^N(\mathbf i\mathcal M^t)_k & = & \sum_{k=0}^N\left(\mathbf i\left(\mathcal A+\mathcal B\right)^t\right)_k\\
& = & (1-r_d)^t(1-r_s)^t\sum_{j=0}^tC_t^j(2r_{div})^{t-j}(1-r_{div})^j\sum_{k=0}^N\left(\mathbf i\mathcal Q_N^{t-j}\right)_k
\end{eqnarray*}
Since $\mathcal Q_N$ is a stochastic matrix, for all $n$, $\mathcal Q_N^n$ is still a stochastic matrix, \emph{i.e.} the entries of each row of $\mathcal Q_N^n$ sum to 1. Therefore:
\begin{eqnarray*}
\sum_{k=0}^N(\mathbf i\mathcal M^t)_k & = & (1-r_d)^t(1-r_s)^t\sum_{j=0}^tC_t^j(2r_{div})^{t-j}(1-r_{div})^j \\
& = & (1-r_d)^t(1-r_s)^t(2r_{div}+1-r_{div})^t = (1-r_d)^t(1-r_s)^t(1+r_{div})^t~,
\end{eqnarray*}
as stated by Equation \eqref{prop:expZt} for $z_0=1$. This result can be easily generalized to the case of $z_0\geq1$ initial B-cells.
\endproof

\section{Proof of Proposition \ref{prop:expselatt}}\label{app:proofprop:expselatt}

\begin{proposition*}
Let us suppose that at time $t=0$ there is a single B-cell entering the GC belonging to the $i^{\textrm{th}}$-affinity class with respect to  the target cell. Moreover, let us suppose that $\mathcal Q_N=R\Lambda_N L$. For all $t\geq1$, the expected number of selected B-cells at time $t$, is:
\begin{equation*}
\E(S_t)=r_s(1-r_s)^{t-1}(1-r_d)^t\displaystyle\sum_{\ell=0}^N(2\lambda_{\ell}r_{div}+1-r_{div})^t \sum_{k=0}^{\overline a_s}r_{i\ell}l_{\ell k}~,
\end{equation*}
\end{proposition*}

\proof
Let us suppose, for the sake of simplicity, that $\mathcal Q_N$ is diagonalizable: 
\begin{equation}\label{eq:Qdiag}
\mathcal Q_N=R\Lambda_N L~,
 \end{equation}
 
We can prove by iteration that:
\begin{equation}\label{eq:Mt}
{\mathcal M}^t=\left(\begin{array}{cc}
\mathcal M_1^t & \displaystyle\sum_{k=0}^{t-1}\mathcal M_1^k\mathcal M_2 \\ \\
\boldsymbol 0_{2\times (N+1)} & \mathcal I_2
\end{array}\right)
\end{equation}

It follows from \eqref{eq:Qdiag} and \eqref{eq:Mt}  that for all $t\geq1$, $\mathcal M^t$ can be written as:
 \begin{equation}\label{eq:Mtdiag}
 \mathcal M^t=\left(\begin{array}{cc}
 RD^tL & \left(R\displaystyle\sum_{k=0}^{t-1}D^kL\right)\mathcal M_2 \\ \\
 \boldsymbol 0_{2\times (N+1)} & \mathcal I_2
 \end{array}\right)~, 
 \end{equation}
 where $D=2(1-r_d)r_{div}(1-r_s)\Lambda_N+(1-r_d)(1-r_{div})(1-r_s)\mathcal I_{N+1}$ is a diagonal matrix. We obtain its expression thanks to Proposition \ref{propM}.\\
 
Moreover, by Proposition \ref{propTildeM} and Equation \eqref{eq:Qdiag} we have:
\begin{equation}\label{eq:tildeMdiag}
\widetilde{\mathcal M}=\left(\begin{array}{cc}
R\widetilde{D}L & \widetilde{\mathcal M}_2 \\
\boldsymbol 0_{2\times (N+1)} & \mathcal I_2
\end{array}\right)~,
\end{equation}
where $\widetilde{D}=2(1-r_d)r_{div}\Lambda_N+(1-r_d)(1-r_{div})\mathcal I_{N+1}$ is a diagonal matrix.\\

Proposition \ref{prop:expMM'} claims:
\begin{equation*}
\E(S_t)=r_s\displaystyle\sum_{k=0}^{\overline a_s}\left(\mathbf i\mathcal M^{t-1}\widetilde{\mathcal M}\right)_k
\end{equation*}
From Equations \eqref{eq:Mtdiag} and \eqref{eq:tildeMdiag}:
\begin{equation*}
\mathcal M^{t-1}\widetilde{\mathcal M}=\left(\begin{array}{ccc}
RD^{t-1}\widetilde{D}L & & RD^{t-1}L\widetilde{\mathcal M}_2+\left(R\displaystyle\sum_{k=0}^{t-2}D^kL\right)\mathcal M_2 \\ \\
\boldsymbol 0_{2\times (N+1)} & & \mathcal I_2
\end{array}\right)
\end{equation*}
Since, by hypothesis, $\vec i=(0,\dots,0,1,0,\dots,0,0)$, with the only 1 being at position $i$, $0\leq i\leq N$, then $\left(\mathbf i\mathcal M^{t-1}\widetilde{\mathcal M}\right)$ denotes the $i^{\textrm{th}}$-row of matrix $\mathcal M^{t-1}\widetilde{\mathcal M}$. Therefore, we are interested in the sum between $0$ and $\overline a_s$ of the elements of the $i^{\textrm{th}}$-row of matrix $\mathcal M^{t-1}\widetilde{\mathcal M}$, \emph{i.e.} of the $i^{\textrm{th}}$-row of matrix $RD^{t-1}\widetilde{D}L$, since clearly $\overline a_s\leq N$.  $D^{t-1}\widetilde{D}$ is a diagonal matrix whose $\ell^{\textrm{th}}$-diagonal element is given by:
\begin{eqnarray*}
\left(D^{t-1}\widetilde{D}\right)_{\ell} & = & (2(1-r_d)r_{div}(1-r_s)\lambda_{\ell}+(1-r_d)(1-r_{div})(1-r_s))^{t-1} \\
& & \cdot(2(1-r_d)r_{div}\lambda_{\ell}+(1-r_d)(1-r_{div})) \\
& = & (1-r_s)^{t-1}(1-r_d)^t\left(2\lambda_{\ell}r_{div}+1-r_{div}\right)^t
\end{eqnarray*}
The result follows observing that: $\left(RD^{t-1}\widetilde{D}L\right)_{ik}=\sum_{\ell=0}^N\left(D^{t-1}\widetilde{D}\right)_{\ell}r_{i\ell}l_{\ell k}$.
\qed
\endproof

\section{Heuristic proof of Proposition \ref{cor:maxrs}}\label{app:proofcor:maxrs}

\begin{proposition*}
For all $t\in\mathbb N$ the value $r_s(t)$ which maximizes the expected number of selected B-cells at the $t^{\textrm{th}}$ maturation cycle is:
\begin{equation*}
r_s(t)=\displaystyle\frac{1}{t}
\end{equation*}
\end{proposition*}

\begin{assumption}\label{hyp1}
$\mathcal Q_N$ converges through its stationary distribution, denoted by $\vec m=(m_i)$, $i\,\in\,\{0,\dots,N\}$.
\end{assumption}
\begin{assumption}\label{hyp2}
$Z_t$ explodes, where $(Z_t)_{t\in\mathbb N}$ is given by Definition \ref{def:Zt}. 
\end{assumption}

Let $\widetilde{Z}_t$, $t\geq0$ be the random variable describing the GC-population size at time $t$ before the selection mechanism is performed for this generation. For the sake of simplicity, let us suppose $\widetilde{Z}_0=1$.
$(\widetilde{Z}_t)_{t\in\mathbb{N}}$ is a MC on $\{0,1,2,\dots\}$. Denoted by $\tilde{p}_k:=\Pro(\widetilde{Z}_1=k)$, $k\,\in\,\{0,1,2\}$:

\begin{equation}
\left\{\begin{array}{l}
\tilde{p}_0=r_d \\
\tilde{p}_1=(1-r_d)(1-r_{div}) \\
\tilde{p}_2=(1-r_d)r_{div}
\end{array}\right.
\end{equation}

It follows: $\tilde{m}:=\E(\widetilde{Z}_1)=(1-r_d)(1-r_{div})+2(1-r_d)r_{div}=(1-r_d)(1+r_{div})$.\\

Conditioning to $Z_t=k$, $\widetilde{Z}_{t+1}$ is distributed as the sum of $k$ independent copies of $\widetilde{Z}_1$, which gives:

\begin{equation}\label{eq:ZtildeSGW}
\E(\widetilde{Z}_t)=\E(Z_{t-1})\E(\widetilde{Z}_1)=\E(Z_1)^{t-1}\E(\widetilde{Z}_1)=(1-r_d)^t(1+r_{div})^t(1-r_s)^{t-1}
\end{equation}

Thanks to Hypotheses \ref{hyp1} and \ref{hyp2}, if $t$ is big enough, 
there is approximately a proportion of $m_i$ elements in the $i^{\textrm{th}}$-affinity class with respect to  $\cible$.
Therefore, on average at time $t$ there are approximately $\sum_{i=0}^{\overline a_s}m_i\E(\widetilde{Z}_t)$ B-cells in the GC belonging to an affinity class with index at most equal to $\overline a_s$ with respect to  $\cible$, before the selection mechanism is performed for this generation. Each one of these cells can be submitted to selection with probability $r_s$, and in this case it will be positively  selected. Hence:

\begin{equation}\label{eq:heurselt}
\E(S_t)\simeq r_s\sum_{i=0}^{\overline a_s}m_i\E(\widetilde{Z}_t)=(1-r_d)^t(1+r_{div})^t(1-r_s)^{t-1}r_s\sum_{i=0}^{\overline a_s}m_i~,
\end{equation}
which is maximized at time $t\geq1$ for $r_s(t)=1/t$.
\begin{remark}
One  observes that the  approximation in  \eqref{eq:heurselt} gives the same value for the optimal $r_s(t)$ as in Proposition \ref{cor:maxrs}. Nevertheless, it does not allow to describe exactly the behavior of $\E(S_t)$, since it is obtained by approximating the distribution of B-cells in the GC with their stationary distribution.
\end{remark}

\end{document}